\documentclass[11pt]{article}
\usepackage{amsmath}
\usepackage{amsthm}
\usepackage{amssymb}
\usepackage{bm}
\usepackage{multirow}

\usepackage{lipsum}
\usepackage{amsfonts}
\usepackage{graphicx}
\usepackage{epstopdf}
\usepackage[symbolgreek,defaultmathsizes]{mathastext}
\usepackage{algorithm}
\usepackage{algorithmic}
\ifpdf
  \DeclareGraphicsExtensions{.eps,.pdf,.png,.jpg}
\else
  \DeclareGraphicsExtensions{.eps}
\fi

\usepackage{mathtools} 
\usepackage{epic}
\usepackage{tikz}
\usetikzlibrary{arrows}

\newtheorem*{remark}{Remark}

\DeclareMathAlphabet{\mathpzc}{OT1}{pzc}{m}{it}

\newcommand{\TheTitle}{Monolithic mixed-dimensional multigrid methods for single-phase flow in fractured porous media}

\title{{\TheTitle}
\footnote{Francisco J. Gaspar has received funding from the European Union's Horizon 2020 research and innovation programme under the Marie Sklodowska-Curie grant agreement No 705402, POROSOS. The work of Andr\'es Arrar\'as, Laura Portero and Carmen Rodrigo is supported in part by the Spanish project FEDER /MCYT MTM2016-75139-R, and the work of Dr. Rodrigo is also supported by the Diputaci\'on General de Arag\'on (Grupo de referencia APEDIF, ref. ).}}

\author{
Andr\'es Arrar\'as\footnote{Departamento de Estad\'{i}stica, Inform\'atica y Matem\'aticas, Universidad P\'ublica de Navarra, Edi\-ficio de Las Encinas, Campus de Arrosad\'ia, 31006 Pamplona, Spain (andres.arraras@unavarra.es).}
\and
  Francisco J. Gaspar\footnote{CWI, Centrum Wiskunde and Informatica, 1098 XG Amsterdam, The Netherlands (gaspar@cwi.nl).}
  \and
   Laura Portero\footnote{Departamento de Estad\'{i}stica, Inform\'atica y Matem\'aticas, Universidad P\'ublica de Navarra, Edi\-ficio de Las Encinas, Campus de Arrosad\'ia, 31006 Pamplona, Spain (laura.portero@unavarra.es).}
  \and 
  Carmen Rodrigo\footnote{IUMA, Departamento de Matem\'atica Aplicada, Universidad de Zaragoza, Pedro Cerbuna 12, 50009 Zaragoza, Spain (carmenr@unizar.es).}
}

\usepackage{amsopn}

\begin{document}

\maketitle

\begin{abstract}
This paper deals with the efficient numerical solution of single-phase flow problems in fractured porous media. A monolithic multigrid method is proposed for solving two-dimensional arbitrary fracture networks with vertical and/or horizontal possibly intersecting fractures. The key point is to combine two-dimensional multigrid components (smoother and inter-grid transfer operators) in the porous matrix with their one-dimensional counterparts within the fractures, giving rise to a mixed-dimensional multigrid method. This combination seems to be optimal since it provides an algorithm whose convergence matches the multigrid convergence factor for solving the Darcy problem. Several numerical experiments are presented to demonstrate the robustness of the monolithic mixed-dimensional multigrid method with respect to the permeability of the fractures, the grid size and the number of fractures in the network.
\end{abstract}



	\section{Introduction}\label{sec:1}
	
	\setcounter{section}{1}
	
The numerical simulation of subsurface flow through fractured porous media is a challenging task which is getting increasing attention in recent years. The essential role played by fractures in different applications --ranging from petroleum extraction to long-term CO$_2$ and nuclear waste storage-- demands the design of efficient discretization methods for solving the corresponding flow models. Depending on the spatial scale under consideration, fractures can be incorporated to such models in essentially two ways. At small scales, when specific locations of micro-fractures are difficult to determine, the so-called dual-porosity models \cite{arb:dou:90,arb:dou:hor:90} are used. In this case, the network of fractures and the bulk or porous matrix are two interacting continua related by a transfer function. On the other hand, at large scales, geological discontinuities represented by localized networks of faults and macro-fractures require the use of discrete fracture models \cite{fle:fum:sco:16,Martin-Jaffre-Roberts}. In these models, fractures can behave either as preferential flow paths or as geological barriers, depending on the permeability contrast between the porous matrix and the fractures themselves.

Henceforth, we shall consider this latter approach. Discrete fracture models typically require fine meshing of the fracture domain to guarantee accurate approximations. Provided that the fracture aperture is small as compared to the characteristic length of the flow domain, this fact can yield computationally expensive discretizations. To avoid such limitations and based on geometrical model reduction techniques, fractures are represented as $(n-1)$-dimensional interfaces immersed into an $n$-dimensional porous matrix. The resulting model is called mixed-dimensional \cite{kei:fum:ber:ste:17,nor:boo:fum:kei:18} or reduced \cite{for:fum:sco:ruf:14,sch:fle:hel:woh:15} model. In this framework, flow can be described by several models within the fractures and in the porous matrix. In \cite{alb:jaf:rob:ser:01,ang:boy:hub:09,boo:nor:yot:18,ang:sco:12,for:fum:sco:ruf:14,Martin-Jaffre-Roberts}, incompressible single-phase Darcy flow is considered in both domains. Extensions to two-phase flow can be found, e.g., in \cite{fum:sco:13,gla:hel:fle:cla:17}. Alternatively, models that consider high-velocity flows within the fractures include Darcy--Brinkman \cite{buk:yot:zun:17,les:ang:qua:11}, Forchheimer \cite{fri:rob:saa:08,kna:rob:14} and Reynolds lubrication \cite{gir:whe:gan:mea:15} equations.

In this paper, we focus on the single-phase Darcy--Darcy coupling between the fractures and the porous matrix. The governing equations comprise a system of mixed-dimensional partial differential equations \cite{boo:nor:vat:17} defined on the $n$-dimensional porous matrix, $(n-1)$-dimensional fractures and $(n-2)$-dimensional intersections between fractures. This problem has been extensively addressed in the literature by means of different discretization techniques. Raviart--Thomas mixed finite element schemes are studied, e.g., in \cite{alb:jaf:rob:ser:01,Martin-Jaffre-Roberts} for the case of conforming meshes on the fracture interfaces. Their extension to non-matching grids is discussed in \cite{fri:mar:rob:saa:12} and, suitably combined with mortar methods, in \cite{boo:nor:yot:18}. The so-called extended finite element methods (XFEM), which permit to mesh the entire domain independently of the fractures, are described in \cite{ang:sco:12,fle:fum:sco:16} and references therein. In addition, further discretization schemes have been proposed for handling general elements and distorted grids, namely: mimetic finite difference methods \cite{ant:for:sco:ver:ver:16,sco:for:sot:17}, discontinuous Galerkin methods \cite{ant:fac:rus:ver:16b}, virtual element methods \cite{ben:ber:pie:sci:14,fum:kei:18}, hybrid high-order methods \cite{cha:pie:for:18}, or multipoint flux approximation methods \cite{ahm:edw:lam:hui:pal:15a,san:ber:nor:12}.

Although a lot of effort has been put into developing numerical schemes for the discretization of fracture models, efficient solvers for the resulting linear systems have not been so deeply investigated. Some relevant works related to this issue include iterative strategies in a domain decomposition framework \cite{ang:sco:12,Martin-Jaffre-Roberts} (see also \cite{fle:fum:sco:16} for a discussion on linear solvers and \cite{hoa:jap:ker:rob:16} for an extension to time-dependent problems), physics-based preconditioners \cite{san:kei:nor:14}, or iterative multiscale methods \cite{haj:kar:jen:2011, ten:saa:haj:2016}. In this context, the aim of this paper is to develop a monolithic multigrid method for solving mixed-dimensional Darcy problems on fractured porous media. To the best of our knowledge, this is the first time that a similar approach is proposed in the literature. For the ease of presentation, we shall assume a distribution of horizontal and vertical fractures that can intersect with each other in virtually any way. This fracture configuration can be efficiently discretized by means of conforming mixed methods based on Raviart--Thomas elements. Further applying suitable quadrature rules, we can derive finite volume schemes that extend the ideas proposed in \cite{rus:whe:83} for non-fractured domains. In addition, we also introduce a novel representation of the network of fractures based on graph theory.

It is well-known that multigrid methods \cite{Bra77, Hackb, Stu_Tro, TOS01, Wess} are among the fastest numerical techniques for solving the large systems of equations arising from the discretization of partial differential equations. They have shown optimal complexity in solving many problems in different areas of application \cite{TOS01}. However, it is the first time that multigrid is applied for solving a mixed-dimensional flow problem in fractured porous media. These algorithms strongly depend on the appropriate choice of their components, mainly the inter-grid transfer operators and the smoother. In this work, a mixed-dimensional multigrid method is proposed to deal with the complex mixed-dimensional problem at once. In a two-dimensional setting, the proposed multigrid solver suitably combines two-dimensional smoother and inter-grid transfer operators in the porous matrix with their one-dimensional counterparts within the fracture network. Due to the saddle point character of the whole resulting system, we choose a multiplicative Schwarz smoother, which has been proved to be efficient for different problems in porous media. The resulting mixed-dimensional monolithic multigrid method shows robustness with respect to the mesh size, the permeability of the fractures, and the number of fractures in the network.

The rest of the paper is organized as follows. Section 2 describes the mixed-dimensional model problem and the spatial discretization considered. In particular, we first focus on the case of a single fracture and then we address the general case of multiple intersecting fractures. Section 3 introduces a mixed-dimensional monolithic multigrid method that combines two-dimensional components in the porous matrix with one-dimensional components in the fractures. Section 4 shows several numerical experiments considering various fracture configurations and permeability distributions that confirm the robustness of the proposed solver. Finally, Section 5 contains some concluding remarks.

	\section{Model problem and discretization}\label{sec:2}
\setcounter{section}{2}

In this section, we introduce the system of equations modeling single-phase Darcy flow in a fractured porous medium. For the ease of presentation, the model is first derived for the case of a single fracture. The weak formulation and its MFE discretization are then obtained. Next, we address the case of intersecting fractures, and emphasize the key points of this extended model. In both cases, the resulting algebraic system for the approximation scheme is provided.

\subsection{A single fracture model problem}

Let $\Omega\subset\mathbb{R}^2$ be an open, bounded, and convex polygonal domain, whose boundary is denoted by $\Gamma=\partial\Omega$. We consider a single-phase, incompressible flow in $\Omega$, governed by Darcy's law in combination with the mass conservation equation, i.e.,
\begin{equation}\label{cont:problem}
\begin{alignedat}{2}
\mathbf{u}&=-\mathbf{K}\nabla p\hspace*{1.5cm}&&\mbox{in }\Omega,\\
\nabla\cdot\mathbf{u}&=q\hspace*{1.5cm}&&\mbox{in }\Omega,\\
p&=0\hspace*{1cm}&&\mbox{on }\Gamma.
\end{alignedat}
\end{equation}
%
Here, $p$ denotes the pressure, $\mathbf{u}$ is the Darcy velocity, $\mathbf{K}\in\mathbb{R}^{2\times 2}$ is the permeability tensor, and $q$ is the source/sink term. We suppose that $\mathbf{K}$ is a diagonal tensor whose entries are strictly positive and bounded in $\Omega$. Homogeneous Dirichlet boundary conditions have been considered for simplicity, but other types of boundary data can also be handled. 

%

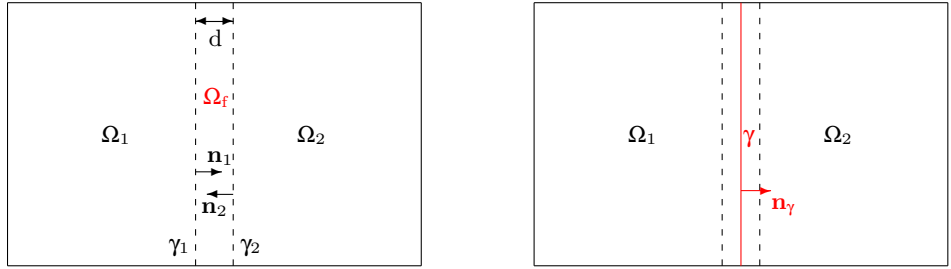
\begin{figure}[t]
	\begin{center}
		\unitlength=0.05cm
		\begin{picture}(260,100)
		\put(5,20){\line(1,0){110}}\put(115,20){\line(0,1){70}}\put(5,90){\line(1,0){110}}\put(5,20){\line(0,1){70}}
		\dashline[20]{2}(55,20)(55,90)
		\dashline[20]{2}(65,20)(65,90)
		\put(30,53){{\footnotesize$\Omega_1$}}
		\put(57,63){{\color{red}\footnotesize$\Omega_f$}} 
		\put(82,53){{\footnotesize$\Omega_2$}}
		\put(55,85){\vector(1,0){10}}
		\put(65,85){\vector(-1,0){10}}
		\put(58.5,78.5){{\footnotesize$d$}}
		\put(48,24){{\footnotesize$\gamma_1$}}
		\put(67,24){{\footnotesize$\gamma_2$}}
		\put(55,45){\vector(1,0){7}}
		\put(58,47.5){\footnotesize$\mathbf{n}_1$}
		\put(65,39){\vector(-1,0){7}}
		\put(56.5,34){\footnotesize$\mathbf{n}_2$}
		
		\put(145,20){\line(1,0){110}}\put(255,20){\line(0,1){70}}\put(145,90){\line(1,0){110}}\put(145,20){\line(0,1){70}}
		\put(200,20){\color{red}\line(0,1){70}}
		\dashline[20]{2}(195,20)(195,90)
		\dashline[20]{2}(205,20)(205,90)
		\put(170,53){{\footnotesize$\Omega_1$}}
		\put(200.7,53){{\color{red}\footnotesize$\gamma$}} 
		\put(222,53){{\footnotesize$\Omega_2$}}
		\put(200,40){\color{red}\vector(1,0){8}}
		\put(208,35){\color{red}\footnotesize$\mathbf{n}_{\gamma}$}
		\end{picture}
	\end{center}\vspace*{-0.7cm}
	\caption{Schematic representation of the original domain (left) and the reduced domain (right).}
	\label{fig:one:fracture}
\end{figure}

To begin with, we suppose that the porous matrix (or bulk) $\Omega$ contains a subset $\Omega_f$ representing a single fracture, which divides the flow domain into two disjoint, connected subdomains $\Omega_1$ and $\Omega_2$, i.e.,
$$
\Omega\backslash\overline{\Omega}_f=\Omega_1\cup\Omega_2,\qquad \Omega_1\cap\Omega_2=\emptyset.
$$
We further denote $\Gamma_{i}=\partial\Omega_i\cap\Gamma$, for $i=1,2,f$, and $\gamma_i=\partial\Omega_i\cap\partial\Omega_f\cap\Omega$, for $i=1,2$. The unit vector normal to $\gamma_i$ pointing outward from $\Omega_i$ is denoted by $\mathbf{n}_i$, for $i=1,2$, see Figure \ref{fig:one:fracture} (left).

If we denote by $p_i$, $\mathbf{u}_i$, $\mathbf{K}_i$ and $q_i$ the restrictions of $p$, $\mathbf{u}$, $\mathbf{K}$ and $q$, respectively, to $\Omega_i$, for $i=1,2,f$, 
then the equations \eqref{cont:problem} are equivalent to the following transmission problem, for $i=1,2,f$ and $j=1,2$,
\begin{subequations}\label{transmission:problem}
	\renewcommand{\theequation}{\theparentequation\alph{equation}}
	\begin{align}
		\mathbf{u}_i&=-\mathbf{K}_i\nabla p_i\hspace*{-1.5cm}&&\hbox{in }\Omega_i,\label{transmission:problem:a}\\
		\nabla\cdot\mathbf{u}_i&=q_i\hspace*{-1.5cm}&&\hbox{in }\Omega_i,\label{transmission:problem:b}\\
		p_j&=p_f\hspace*{-1.5cm}&&\hbox{on }\gamma_j,\label{transmission:problem:c}\\
		\mathbf{u}_j\cdot\mathbf{n}_j&=\mathbf{u}_f\cdot\mathbf{n}_j\hspace*{-1.5cm}&&\hbox{on }\gamma_j,\label{transmission:problem:d}\\
		p_i&=0\hspace*{-1.5cm}&&\hbox{on }\Gamma_{i}.\label{transmission:problem:e}
	\end{align}
\end{subequations}
Note that the equations \eqref{transmission:problem:c} and \eqref{transmission:problem:d} provide coupling conditions that guarantee the continuity of the pressure and the normal flux, respectively, across the interfaces between $\Omega_f$ and $\Omega_i$, for $i=1,2$.

The model provided by the equations \eqref{transmission:problem:a}-\eqref{transmission:problem:e} is sometimes referred to as equi-dimensional model \cite{fle:fum:sco:16}, and assumes that both the bulk and the fracture domains share the same dimension. As an alternative to this model, we shall define the so-called mixed-dimensional or reduced model, in which the fracture is viewed as a manifold of co-dimension one (that is, an interface between the bulk subdomains $\Omega_1$ and $\Omega_2$). Based on model reduction techniques, this idea was first proposed in \cite{Martin-Jaffre-Roberts} and is commonly used in the framework of fractured porous media \cite{ami:ker:mar:rob:2006,delpra:fum:sco:17,fle:fum:sco:16,for:fum:sco:ruf:14,fri:mar:rob:saa:12}. From a numerical viewpoint, the mixed-dimensional approach avoids fine meshing of the fracture domain, thus reducing the computational cost of the overall discretization.


According to \cite{Martin-Jaffre-Roberts}, there exists a non-self-intersecting one-dimensional manifold $\gamma$ such that the fracture can be expressed as
$$
\Omega_f=\left\{\mathbf{x}\in\Omega:\mathbf{x}=\mathbf{s}+\theta\,\mathbf{n}_{\gamma}, \hbox{ for some } \mathbf{s}\in\gamma \hbox{ and } |\theta|<\textstyle\dfrac{d(\mathbf{s})}{2}\right\},
$$
where $d(\mathbf{s})>0$ denotes the thickness of the fracture at $\mathbf{s}$ in the normal direction, and $\mathbf{n}_{\gamma}$ is the outward unit normal to $\gamma$ with a fixed orientation from $\Omega_1$ to $\Omega_2$. Note that, with this definition, $\mathbf{n}_{\gamma}=\mathbf{n}_1=-\mathbf{n}_2$ (see Figure \ref{fig:one:fracture}). We will assume that the thickness is smaller than the other characteristic dimensions of the fracture.

The key point in this procedure is to collapse the fracture $\Omega_f$ into the line $\gamma$, and integrate the equations \eqref{transmission:problem:a} and \eqref{transmission:problem:b} for the index $f$ along the fracture thickness. In doing so, we need to split up such equations into their normal and tangential parts. Let us denote the projection operators onto the normal and tangent spaces of $\gamma$ as $\mathbf{P}_{\mathbf{n}}=\mathbf{n}_{\gamma}\mathbf{n}_{\gamma}^T$ and $\mathbf{P}_{\boldsymbol{\tau}}=\mathbf{I}-\mathbf{P}_{\mathbf{n}}$, $\mathbf{I}$ being the identity tensor. For regular vector- and scalar-valued functions $\mathbf{g}$ and $g$, the tangential divergence and gradient operators on the fracture are defined, respectively, as
\begin{equation}\label{tangential:div:grad}
\nabla^{\boldsymbol{\tau}}\cdot\mathbf{g}=\mathbf{P}_{\boldsymbol{\tau}}:\mathbf{\nabla}\mathbf{g},\qquad\nabla^{\boldsymbol{\tau}} g=\mathbf{P}_{\boldsymbol{\tau}}\nabla g.
\end{equation}
Following \cite{for:fum:sco:ruf:14}, we assume that the permeability tensor $\mathbf{K}_f$ decomposes additively as
\begin{equation}\label{K:decomposition}
\mathbf{K}_f=K_{f}^{\mathbf{n}}\mathbf{P}_{\mathbf{n}}+K_{f}^{\boldsymbol{\tau}}\mathbf{P}_{\boldsymbol{\tau}},
\end{equation}
where $K_{f}^{\mathbf{n}}$ and $K_{f}^{\boldsymbol{\tau}}$ are defined to be strictly positive and bounded in $\Omega_f$.

In this framework, we introduce the so-called reduced variables, namely: the reduced pressure $p_{\gamma}$, the reduced Darcy velocity $\mathbf{u}_{\gamma}$, and the reduced source/sink term $q_{\gamma}$, formally defined as \cite{delpra:fum:sco:17,Martin-Jaffre-Roberts}
$$
p_{\gamma}(\mathbf{s})=\textstyle\dfrac{1}{d(\mathbf{s})}(p_f,1)_{\ell(\mathbf{s})},\quad\
\mathbf{u}_{\gamma}(\mathbf{s})=(\mathbf{P}_{\boldsymbol{\tau}}\mathbf{u}_f,1)_{\ell(\mathbf{s})},\quad\
q_{\gamma}(\mathbf{s})=(q_f,1)_{\ell(\mathbf{s})},
$$
where $\ell(\mathbf{s})=\left(-\frac{d(\mathbf{s})}{2},\frac{d(\mathbf{s})}{2}\right)$. 
Hence, we obtain the following interface problem, for $i=1,2$,
\begin{subequations}
	\renewcommand{\theequation}{\theparentequation\alph{equation}}\label{interface:problem}
	\begin{align}
		\mathbf{u}_i&=-\mathbf{K}_i\nabla p_i&&\hbox{in }\Omega_i,\label{interface:problem:a}\\
		\nabla\cdot\mathbf{u}_i&=q_i&&\hbox{in }\Omega_i,\label{interface:problem:b}\\
		\mathbf{u}_{\gamma}&=-dK_{f}^{\boldsymbol{\tau}}\nabla^{\boldsymbol{\tau}} p_{\gamma}&&\hbox{on }\gamma,\label{interface:problem:c}\\
		\nabla^{\boldsymbol{\tau}}\cdot\mathbf{u}_{\gamma}&=q_{\gamma}+(\mathbf{u}_1\cdot\mathbf{n}_1+\mathbf{u}_2\cdot\mathbf{n}_2)&&\hbox{on }\gamma,\label{interface:problem:d}\\
		\alpha_{\gamma}(p_i-p_{\gamma})&=\xi\,\mathbf{u}_i\cdot\mathbf{n}_{i}-(1-\xi)\,\mathbf{u}_{i+1}\cdot\mathbf{n}_{i+1}&&\hbox{on }\gamma,\label{interface:problem:e}\\
		p_i&=0&&\hbox{on }\Gamma_i,\label{interface:problem:f}\\
		p_{\gamma}&=0&&\hbox{on }\partial\gamma,\label{interface:problem:g}
	\end{align}
\end{subequations}
where $\alpha_{\gamma}=2K_{f}^{\mathbf{n}}/d$ and the index $i$ is supposed to vary in $\mathbb{Z}/2\mathbb{Z}$, so that, if $i=2$, then $i+1=1$. Following \cite{ang:boy:hub:09,Martin-Jaffre-Roberts}, $\xi\in(1/2,1]$ is a closure parameter related to the pressure cross profile in the fracture. The ratio $K_{f}^{\mathbf{n}}/d$ and the product $K_{f}^{\boldsymbol{\tau}}d$ are sometimes referred to as effective permeabilities in the normal and tangential directions to the fracture, respectively \cite{fle:fum:sco:16}.

In the preceding system, \eqref{interface:problem:c} represents Darcy's law in the tangential direction of the fracture, while \eqref{interface:problem:d} models mass conservation inside the fracture. Remarkably, the additional source term $\mathbf{u}_1\cdot\mathbf{n}_1+\mathbf{u}_2\cdot\mathbf{n}_2$ is introduced on $\gamma$ to take into account the contribution of the subdomain flows to the fracture flow. In turn, \eqref{interface:problem:e} is obtained by averaging Darcy's law in the normal direction to the fracture and using a quadrature rule with weights $\xi$ and $1-\xi$ for integrating $\mathbf{u}_f\cdot\mathbf{n}_i$ across the fracture, for $i=1,2$. Formally, it can be regarded as a Robin boundary condition for the subdomain $\Omega_i$ that involves the pressure in the fracture $p_{\gamma}$ and the normal flux from the neighboring subdomain $\Omega_{i+1}$. It is quite usual to express \eqref{interface:problem:e} in terms of average operators for the pressures and normal fluxes, and jump operators for the pressures across the fracture \cite{ang:sco:12,delpra:fum:sco:17}. 


\subsection{Weak formulation} In this subsection, we present the weak formulation of the interface problem stated above. To this end, we first introduce the following function spaces 
\begin{align*}
	%
	%
	\mathbf{W}&=\{\mathbf{v}=(\mathbf{v}_1,\mathbf{v}_2,\mathbf{v}_{\gamma})\in H(\mathrm{div},\Omega_1)\times H(\mathrm{div},\Omega_2)\times H(\mathrm{div}^{\boldsymbol{\tau}},\gamma):\mathbf{v}_i\cdot\mathbf{n}_{i}\in L^2(\gamma),
	\mbox{ for } i=1,2\},\\[1ex]
	M&=\{r=(r_1,r_2,r_{\gamma})\in L^2(\Omega_1)\times L^2(\Omega_2)\times L^2(\gamma)\},	
\end{align*}
endowed with the norms \cite{Martin-Jaffre-Roberts}
\begin{align*}
	\|\mathbf{v}\|^2_{\mathbf{W}}&=\sum_{i=1}^2\left(\|\mathbf{v}_i\|^2_{L^2(\Omega_i)}+\|\nabla\cdot\mathbf{v}_i\|^2_{L^2(\Omega_i)}+\|\mathbf{v}_i\cdot\mathbf{n}_{i}\|^2_{L^2(\gamma)}\right)
	+\|\mathbf{v}_{\gamma}\|^2_{L^2(\gamma)}+\|\nabla^{\boldsymbol{\tau}}\cdot\mathbf{v}_{\gamma}\|^2_{L^2(\gamma)},\\
	\|r\|^2_M&=\sum_{i=1}^2\|r_i\|^2_{L^2(\Omega_i)}+\|r_{\gamma}\|^2_{{L^2(\gamma)}}.
\end{align*}
%
Here, we use the well-known spaces
\begin{align*}
	H(\mbox{div},\Omega_i)&=\{\mathbf{v}_i\in(L^2(\Omega_i))^2:\nabla\cdot\mathbf{v}_i\in L^2(\Omega_i)\},\qquad i=1,2,\\
	H(\mathrm{div}^{\boldsymbol{\tau}},\gamma)&=\{\mathbf{v}_{\gamma}\in (L^2(\gamma))^2:\nabla^{\boldsymbol{\tau}}\cdot\mathbf{v}_{\gamma}\in L^2(\gamma)\},
\end{align*}
and assume that the elements $\mathbf{v}_{\gamma}\in H(\mathrm{div}^{\boldsymbol{\tau}},\gamma)$ are aligned with $\gamma$, that is, $\mathbf{v}_{\gamma}=v_{\gamma}\boldsymbol{\tau}$, where $\boldsymbol{\tau}$ denotes the tangent vector to $\gamma$.	Note that, in order to take into proper account the Robin boundary condition, we need more regularity in $\mathbf{W}$ than the usual $H(\mathrm{div},\cdot)$-regularity commonly used for weak formulations in the context of mixed finite element methods \cite{Martin-Jaffre-Roberts}.

Let $a:\mathbf{W}\times\mathbf{W}\rightarrow\mathbb{R}$ and $b:\mathbf{W}\times M\rightarrow\mathbb{R}$ be the bilinear forms defined by
\begin{align}
	a(\mathbf{u},\mathbf{v})&=\sum_{i=1}^2\left(\mathbf{K}_i^{-1}\mathbf{u}_i,\mathbf{v}_i\right)_{\Omega_i}+\left((dK_{f}^{\boldsymbol{\tau}})^{-1}\mathbf{u}_{\gamma},\mathbf{v}_{\gamma}\right)_{\gamma}
	+\sum_{i=1}^2\left(\alpha_{\gamma}^{-1}(\xi\,\mathbf{u}_i\cdot\mathbf{n}_{i}-(1-\xi)\,\mathbf{u}_{i+1}\cdot\mathbf{n}_{i+1}),\mathbf{v}_i\cdot\mathbf{n}_{i}\right)_{\gamma},\nonumber\\
	b(\mathbf{u},r)&=\sum_{i=1}^2\left(\nabla\cdot \mathbf{u}_i,r_i\right)_{\Omega_i}+\left(\nabla^{\tau}\cdot \mathbf{u}_{\gamma},r_{\gamma}\right)_{\gamma}
	-\left(\mathbf{u}_1\cdot\mathbf{n}_1+\mathbf{u}_2\cdot\mathbf{n}_2,r_{\gamma}\right)_{\gamma}.\nonumber
\end{align}
Accordingly, let $L:M\rightarrow\mathbb{R}$ be the linear form associated with the source terms, i.e.,
$$
L(r)=\sum_{i=1}^2\left(q_i,r_i\right)_{\Omega_i}+\left(q_{\gamma},r_{\gamma}\right)_{\gamma}.
$$
In this framework, the weak formulation of the interface problem  \eqref{interface:problem} reads: \emph{Find $(\mathbf{u},p)\in\mathbf{W}\times M$ such that}
\begin{equation}\label{weak:form}
\begin{alignedat}{2}
a(\mathbf{u},\mathbf{v})-b(\mathbf{v},p)&=0 &&\hspace*{1cm}\forall\,\mathbf{v}\in\mathbf{W},\\
b(\mathbf{u},r)&=L(r)&&\hspace*{1cm}\forall\,r\in M.
\end{alignedat}
\end{equation}
The existence and uniqueness of solution to this problem is proved in \cite{Martin-Jaffre-Roberts} for the case $\xi>1/2$, assuming that the permeabilities in both subdomains and the coefficients $K_{f}^{\mathbf{n}}/d$ and $K_{f}^{\boldsymbol{\tau}}d$ are bounded by positive constants.


\subsection{Mixed finite element approximation}\label{subsec:mfe}

Let us assume that the subdomains $\Omega_i$ admit rectangular partitions $\mathcal{T}^{h}_{i}$, for $i=1,2$, that match at the interface $\gamma$. Such meshes $\mathcal{T}^{h}_{i}$ induce a unique partition on $\gamma$ denoted by $\mathcal{T}^h_{\gamma}$.

Let $\mathbf{W}^h_{i}\times M^h_{i}$ be the lowest order Raviart--Thomas mixed finite element spaces defined on $\mathcal{T}^h_{i}$, for $i=1,2,\gamma$, and let us introduce the global spaces
\begin{equation*}
	\mathbf{W}^h=\displaystyle\bigoplus_{i=1,2,\gamma}\mathbf{W}^{h}_{i},\qquad\qquad M^h=\displaystyle\bigoplus_{i=1,2,\gamma} M^{h}_{i}.
\end{equation*}
Following \cite{arb:whe:yot:97,rus:whe:83}, we will use numerical quadrature rules for evaluating some of the integrals in \eqref{weak:form}. More specifically, based on the expression for $a(\mathbf{u},\mathbf{v})$ defined above, we set the following discrete bilinear form
\begin{align*}
	a_{h}(\mathbf{u},\mathbf{v})=&\sum_{i=1}^2\left(\mathbf{K}_i^{-1}\mathbf{u}_i,\mathbf{v}_i\right)_{\Omega_i,\mathbf{TM}}+\left((d \,K_{f}^{\boldsymbol{\tau}})^{-1}\mathbf{u}_{\gamma},\mathbf{v}_{\gamma}\right)_{\gamma,\mathbf{T}}\\
	&\hspace*{2cm}+\sum_{i=1}^2\left(\alpha_{\gamma}^{-1}(\xi\,\mathbf{u}_i\cdot\mathbf{n}_i-(1-\xi)\,\mathbf{u}_{i+1}\cdot\mathbf{n}_{i+1}),\mathbf{v}_i\cdot\mathbf{n}_i\right)_{\gamma},
\end{align*}
where $(\cdot,\cdot)_{\gamma,\mathbf{T}}$ denotes the application of the trapezoidal rule for computing the inner-product integral over $\gamma$, and $(\cdot,\cdot)_{\Omega_i,\mathbf{TM}}$ is defined, for any vector-valued functions $\mathbf{f},\mathbf{g}\in\mathbb{R}^2$, as \cite{arb:whe:yot:97,rus:whe:83}
\begin{equation*}
	(\mathbf{f},\mathbf{g})_{\Omega_i,\mathbf{TM}}=(f_1,g_1)_{\Omega_i,\mathbf{T}\times\mathbf{M}}+(f_2,g_2)_{\Omega_i,\mathbf{M}\times \mathbf{T}}.
\end{equation*}
In this case, the integral of the $i$th component of the vectors, for $i=1,2$, is computed by using the trapezoidal rule in the $i$th direction and the midpoint rule in the other direction. On the other hand, the discrete counterpart to $L(r)$ is given by
\begin{equation*}
	L_h(r)=\sum_{i=1}^2\left(q_i,r_i\right)_{\Omega_i,\mathbf{M}}+\left(q_{\gamma},r_{\gamma}\right)_{\gamma,\mathbf{M}},
\end{equation*}
where $(\cdot,\cdot)_{G,\mathbf{M}}$ means the application of the midpoint rule for computing the corresponding inner-product integral over $G$. In this context, the mixed finite element approximation to \eqref{weak:form} may be written as: \emph{Find $(\mathbf{u}_h,p_h)\in\mathbf{W}^h\times M^h$ such that}
\begin{equation}\label{mfe:approx}
\begin{alignedat}{2}
a_{h}(\mathbf{u}_h,\mathbf{v}_h)-b(\mathbf{v}_h,p_h)&=0&&\hspace*{1cm}\forall\,\mathbf{v}_h\in\mathbf{W}^h,\\
b(\mathbf{u}_h,r_h)&=L_h(r_h)&&\hspace*{1cm}\forall\,r_h\in M^h.
\end{alignedat}
\end{equation}
Note that the definition of the global spaces $\mathbf{W}^h$ and $M^h$ implies $\mathbf{u}_h=(\mathbf{u}_1^h,\mathbf{u}_2^h,\mathbf{u}_{\gamma}^h)$ and $p_h=(p_1^h,p_2^h,p_{\gamma}^h)$. Following \cite{rus:whe:83}, it can be proved that this method is closely related to the so-called two-point flux approximation (TPFA) method \cite{eym:gal:gui:her:mas:14}. 

\subsection{Algebraic linear system} Let $\{\mathbf{v}_{k,i}\}_{i=1}^{E_k}$ and $\{r_{k,i}\}_{i=1}^{C_k}$ denote the basis functions of $\mathbf{W}_{k}^h$ and $M_{k}^h$, respectively, for $k=1,2,\gamma$. Here, $E_k$ and $C_k$ stand for the number of edges and cells in $\mathcal{T}_{k}^h$, respectively.	Thus, the unknowns in \eqref{mfe:approx} can be expressed as
$$
\mathbf{u}^h_{k}=\sum_{i=1}^{E_k}U_{k,i}\,\mathbf{v}_{k,i},\qquad p^h_{k}=\sum_{i=1}^{C_k}P_{k,i}\,r_{k,i},
$$
for $k=1,2,\gamma$. If we define the vectors $U_k\in\mathbb{R}^{E_k}$ and $P_k\in\mathbb{R}^{C_k}$ with components $U_{k,i}$ and $P_{k,i}$, respectively, for $k=1,2,\gamma$, then the algebraic linear system stemming from \eqref{mfe:approx} is a saddle-point problem of the form
$$\begin{bmatrix}
A_1 & D^T & 0 & B_1^T & 0 & F_1^T\\[0.5ex]
D & A_2 & 0 & 0 & B_2^T & F_2^T\\[0.5ex]
0 & 0 & A_{\gamma} & 0 & 0 & B_{\gamma}^T\\[0.5ex]
B_1 & 0 & 0 & 0 & 0 & 0\\[0.5ex]
0 & B_2 & 0 & 0 & 0 & 0\\[0.5ex]
F_1 & F_2 & B_{\gamma} & 0 & 0 & 0
\end{bmatrix}
\begin{bmatrix}
U_1\\[0.5ex]U_2\\[0.5ex]U_{\gamma}\\[0.5ex]P_1\\[0.5ex]P_2\\[0.5ex]P_{\gamma}
\end{bmatrix}=
\begin{bmatrix}
0\\[0.5ex]0\\[0.5ex]0\\[0.5ex]Q_1\\[0.5ex]Q_2\\[0.5ex]Q_{\gamma}
\end{bmatrix}.
$$
In particular, the entries of the matrices $A_k\in\mathbb{R}^{E_k\times E_k}$, $B_k\in\mathbb{R}^{C_k\times E_k}$ and $F_{k}\in\mathbb{R}^{C_{\gamma}\times E_k}$, for $k=1,2$, $D\in\mathbb{R}^{E_2\times E_1}$, $A_\gamma\in\mathbb{R}^{E_{\gamma}\times E_{\gamma}}$ and $B_\gamma\in\mathbb{R}^{C_{\gamma}\times E_{\gamma}}$ are
\begin{align*}
	[A_k]_{i,j}&=\left(\mathbf{K}_k^{-1}\mathbf{v}_{k,j},\mathbf{v}_{k,i}\right)_{\Omega_k,\mathbf{TM}}+\alpha_{\gamma}^{-1}\left(\xi\mathbf{v}_{k,j}\cdot\mathbf{n}_k,\mathbf{v}_{k,i}\cdot\mathbf{n}_k\right)_{\gamma},&& k=1,2,\\[0.5ex]
	[B_k]_{i,j}&=-(r_{k,i},\nabla\cdot\mathbf{v}_{k,j})_{\Omega_k},&& k=1,2,\\[0.5ex]
	[F_k]_{i,j}&=\left(r_{\gamma,i},\mathbf{v}_{k,j}\cdot\mathbf{n}_k\right)_{\gamma},&& k=1,2,\\[0.5ex]
	[D]_{i,j}&=\,\alpha_{\gamma}^{-1}\left((\xi-1)\mathbf{v}_{2,i}\cdot\mathbf{n}_2,\mathbf{v}_{1,j}\cdot\mathbf{n}_1\right)_{\gamma},\\[0.5ex]
	[A_{\gamma}]_{i,j}&=\left((dK_{f}^{\boldsymbol{\tau}})^{-1}\mathbf{v}_{\gamma,j},\mathbf{v}_{\gamma,i}\right)_{\gamma,\mathbf{T}},\\[0.5ex]
	[B_{\gamma}]_{i,j}&=-\left(r_{\gamma,i},\nabla^{\tau}\cdot\mathbf{v}_{\gamma,j}\right)_{\gamma},
\end{align*}
where $[\cdot]_{i,j}$ indicates the $(i,j)$th element of the matrix. Note that the use of the quadrature rules $(\cdot,\cdot)_{\Omega_k,\mathbf{TM}}$ and $(\cdot,\cdot)_{\gamma,\mathbf{T}}$ diagonalizes the matrices $A_1$, $A_2$ and $A_{\gamma}$. In turn, the components of the vectors $Q_k\in\mathbb{R}^{C_k}$, for $k=1,2,\gamma$, are given by
\begin{align*}
	[Q_k]_j&=-(q_k,r_{k,j})_{\Omega_k,\mathbf{M}},&&\hspace*{-2cm} k=1,2,\\[0ex]
	[Q_{\gamma}]_j&=-(q_{\gamma},r_{\gamma,j})_{\gamma,\mathbf{M}},
\end{align*}
where $[\cdot]_{j}$ denotes the $j$th component of the vector.

Next subsection is devoted to the general case of multiple intersecting fractures. With the aim of achieving a unified notation, we shall group the unknowns corresponding to the two-dimensional bulk subdomains into the vectors $U^2=[U_1,U_2]^T$ and $P^2=[P_1,P_2]^T$, and those associated with the one-dimensional fractures into $U^1=U_{\gamma}$ and $P^1=P_{\gamma}$. In these cases, the superscripts stand for the corresponding dimensions. Using such notations, the preceding system can be rewritten as
\begin{equation}
\label{linear:system:1fracture}
\begin{bmatrix}
A_{2,2} &  0 & B_{2,2}^T &  F_{2,1}^T\\[0.5ex]
0 &  A_{1,1} & 0 &  B_{1,1}^T\\[0.5ex]
B_{2,2} & 0 & 0 & 0 \\[0.5ex]
F_{2,1} &  B_{1,1} & 0 & 0 
\end{bmatrix}
\begin{bmatrix}
U^2\\[0.5ex]U^1\\[0.5ex]P^2\\[0.5ex]P^1
\end{bmatrix}=
\begin{bmatrix}
0\\[0.5ex]0\\[0.5ex]Q^2\\[0.5ex]Q^1
\end{bmatrix},
\end{equation}
where, accordingly, $Q^2=[Q_1,Q_2]^T$ and $Q^1=Q_{\gamma}$.

\subsection{The case of multiple intersecting fractures}

Let us now consider problem \eqref{cont:problem} posed on a geological domain
subdivided into $m$ subdomains $\Omega_i$, for $i\in \mathcal{I}_2=\{ 1, 2, \ldots, m\}$, naturally separated by a collection of fractures $\Omega_{i,j}$, for $(i,j)\in\mathcal{I}_1$. Here and henceforth, $\mathcal{I}_1$ is a set of indices $(i,j)$, with $i$, $j\in\mathcal{I}_2$ and $i<j$, that permits us to label the fractures. In particular, $\Omega_{i,j}$ denotes a fracture that is adjacent to subdomains $\Omega_i$ and $\Omega_j$. In this framework, it holds
$$\Omega\,\backslash \left(\displaystyle\bigcup_{(i,j)\in \mathcal{I}_1}\overline{\Omega}_{i,j}\right)=\displaystyle\bigcup_{i\in \mathcal{I}_2}\Omega_i,$$
with $\Omega_i\cap\Omega_j=\emptyset$, for $i\neq j$.

Following the ideas of the preceding subsection, let us suppose that there exist certain non-self-intersecting one-dimensional manifolds $\gamma_{i,j}$ such that the fractures can be defined as
$$\Omega_{i,j}=\left\{\mathbf{x}\in\Omega:\mathbf{x}=\mathbf{s}+\theta\,\mathbf{n}_{i,j}, \hbox{ for some } \mathbf{s}\in\gamma_{i,j} \hbox{ and } |\theta| <\frac{d_{i,j}(\mathbf{s})}{2}\right\},$$
where $d_{i,j}(\mathbf{s})$ denotes the thickness of the fracture $\Omega_{i,j}$ at $\mathbf{s}$ in the normal direction and $\mathbf{n}_{i,j}$ is the outward unit normal to $\gamma_{i,j}$ oriented from $\Omega_i$ to $\Omega_{j}$, for $(i,j)\in \mathcal{I}_1$. Let us denote by $\{\sigma_{i,j,k}\}_{(i,j,k)\in \mathcal{I}_0^T}$ the $T$-shaped intersections of three fractures (i.e., $\sigma_{i,j,k}=\overline{\gamma}_{i,j}\cap\overline{\gamma}_{j,k}\cap\overline{\gamma}_{i,k}$), and by  $\{\sigma_{i,j,k,l}\}_{(i,j,k,l)\in \mathcal{I}_0^X}$ the $X$-shaped intersections of four fractures (for instance, $\sigma_{i,j,k,l}=\overline{\gamma}_{i,j}\cap\overline{\gamma}_{j,k}\cap\overline{\gamma}_{k,l}\cap\overline{\gamma}_{i,l}$).

%

\begin{figure}[t]
	\begin{center}
		\begin{minipage}{0.4\textwidth}
			\begin{center}
				\unitlength=0.05cm
				\begin{picture}(100,100)
				\put(0,0){\line(1,0){100}}\put(0,0){\line(0,1){100}}\put(100,0){\line(0,1){100}}
				\put(0,100){\line(1,0){100}}
				\put(50,0){\line(0,1){100}}
				\put(0,50){\line(1,0){100}}
				\put(50,75){\line(1,0){50}}
				\put(75,50){\line(0,1){50}}
				\put(50,62.5){\line(1,0){25}}
				\put(62.5,50){\line(0,1){25}}
				\put(48.4,61.05){\footnotesize$\bullet$}
				\put(48.4,73.55){\footnotesize$\bullet$}
				\put(61.12,73.55){\footnotesize$\bullet$}
				\put(61.1,48.5){\footnotesize$\bullet$}
				\put(73.5,48.5){\footnotesize$\bullet$}
				\put(73.5,61.05){\footnotesize$\bullet$}
				\put(47.8,48.6){\footnotesize$\times$}
				\put(60.3,61.1){\footnotesize$\times$}
				\put(72.9,73.6){\footnotesize$\times$}
				\put(23,73){\footnotesize$\Omega_{1}$}
				\put(59,85){\footnotesize$\Omega_{2}$}
				\put(84,85){\footnotesize$\Omega_{3}$}
				\put(84,60){\footnotesize$\Omega_{4}$}
				\put(53,54.5){\footnotesize$\Omega_{6}$}
				\put(65.5,54.5){\footnotesize$\Omega_{5}$}
				\put(53,67){\footnotesize$\Omega_{7}$}
				\put(65.5,67){\footnotesize$\Omega_{8}$}
				\put(73,23){\footnotesize$\Omega_{9}$}
				\put(22,23){\footnotesize$\Omega_{10}$}
				\end{picture}
			\end{center}
		\end{minipage}\hspace*{1.5cm}
		\begin{minipage}{0.4\textwidth}
			
			\begin{center}
				\begin{tikzpicture}[scale=2.1]
				\tikzstyle{every node}=[draw,shape=circle];
				\path (180:1cm) node (1) {$1$};
				\path (216:1cm) node (10) {$10$};
				\path (252:1cm) node (9) {$9 $};
				\path (288:1cm) node (8) {$8 $};
				\path (324:1cm) node (7) {$7 $};
				\path (0:1cm) node (6) {$6 $};
				\path (36:1cm) node (5) {$ 5$};
				\path (72:1cm) node (4) {$4 $};
				\path (108:1cm) node (3) {$3 $};
				\path (144:1cm) node (2) {$2 $};
				\draw (1) -- (2)
				(1) -- (6)
				(1) -- (7)
				(1) -- (10)
				(2) -- (3)
				(2) -- (7)
				(2) -- (8)
				(3) -- (4)
				(4) -- (5)
				(4) -- (8)
				(4) -- (9)
				(5) -- (6)
				(5) -- (8)
				(5) -- (9)
				(6) -- (7)
				(6) -- (9)
				(7) -- (8)
				(9) -- (10);
				
				\end{tikzpicture}
			\end{center}
		\end{minipage}
	\end{center}\vspace*{0.3cm}
	\caption{Schematic representation of the reduced domain (left) and the associated graph (right).}
	\label{fig:graph}
\end{figure}
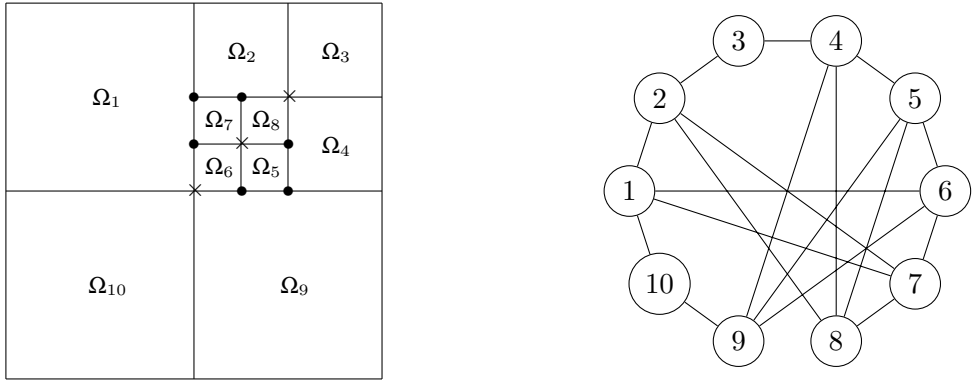

In the sequel, we present an example to illustrate the used notations for a flow domain containing both horizontal and vertical fractures that may intersect. Figure \ref{fig:graph} (left) shows the schematic representation of a domain that contains $m=10$ subdomains and 18 fractures. This geometry will be later considered as a benchmark problem in the section devoted to the numerical experiments (cf. Subsection \ref{subsec:benchmark}).

In general, given $i=1,2,\ldots,m-1$, let us define the set $\mathcal{N}_i$ that contains the indices of the subdomains $\Omega_j$, with $j>i$, that are adjacent to $\Omega_i$. Specifically, in this example, we have
$$\begin{array}{lll}
\mathcal{N}_1=\{2,6,7,10\},\quad&
\mathcal{N}_4=\{5,8,9\},\quad&
\mathcal{N}_7=\{8\},\\
\mathcal{N}_2=\{3,7,8\},&
\mathcal{N}_5=\{6,8,9\},&
\mathcal{N}_8=\emptyset,\\
\mathcal{N}_3=\{4\},&
\mathcal{N}_6=\{7,9\},&
\mathcal{N}_9=\{10\}.
\end{array}$$
Then, the set of indices $\mathcal{I}_1$ denoting the one-dimensional collapsed fractures may be defined as
$$
\mathcal{I}_1=\{(i,j)\in\mathbb{N}^2: i\in\{1,2,\ldots,m-1\},\,j\in\mathcal{N}_i\}.
$$
It is straightforward to see that the number of fractures, in this case 18, is equal to $\sum_{i=1}^{m-1}|\mathcal{N}_i|$, where $|\mathcal{N}_i|$ denotes the cardinal of the set $\mathcal{N}_i$. In this case, $\mathcal{I}_1$ is given by
\begin{align*}\mathcal{I}_1=\{&(1,2),(1,6),(1,7),(1,10),(2,3),(2,7),(2,8),(3,4),(4,5),(4,8),(4,9),\\&
	(5,6),(5,8),(5,9),(6,7),(6,9),(7,8),(9,10)\}.
\end{align*}

In order to define the sets of indices $\mathcal{I}_0^T$ and $\mathcal{I}_0^X$ that refer to the zero-dimensional intersections of fractures, we will introduce a suitable graph representing the problem. In particular, the graph assigns a node to each subdomain and considers a path connecting the nodes $i$ and $j$, as long as $j\in\mathcal{N}_i$ (or, equivalently, whenever the pair of indices $(i,j)\in \mathcal{I}_1$). In other words, the graph nodes stand for the subdomains $\{\Omega_i\}_{i=1}^{m}$, and the paths represent the collapsed fractures $\{\gamma_{i,j}\}_{(i,j)\in \mathcal{I}_1}$. Figure \ref{fig:graph} (right) shows the graph corresponding to the reduced domain on the left. In this framework, the set of indices $\mathcal{I}_0^T$ associated to $T$-shaped intersections of three fractures is defined as
\begin{align*}\mathcal{I}_0^T=\{&(i,j,k)\in\mathbb{N}^3: i,j,k\in \mathcal{I}_2, \hbox{ with } i<j<k, \hbox{ such that the graph contains}\\
&\hspace*{0cm}\hbox{a closed path passing through the nodes } i,j,k\}.
\end{align*}
In the example under consideration, the $T$-shaped intersections --marked with a bullet in Figure \ref{fig:graph} (left)-- are given by the set of indices
$$\mathcal{I}_0^T=\{(1,2,7),(1,6,7),(2,7,8),(4,5,8),(4,5,9),(5,6,9)\}.$$
On the other hand, the set of indices $\mathcal{I}_0^X$ corresponding to $X$-shaped intersections of four fractures is defined as
\begin{align*}
\mathcal{I}_0^X=\{&(i,j,k,l)\in\mathbb{N}^4: i,j,k,l\in \mathcal{I}_2, \hbox{ with } i<j<k<l, \hbox{ such that the graph}\\
&\hspace*{0cm}\hbox{contains a closed path passing through the nodes } i,j,k,l, \hbox{ and } (i,j,k),\\
&\hspace*{0cm}(i,j,l),(i,k,l),(j,k,l)\notin\mathcal{I}_0^T\}.
\end{align*}
In the example, the $X$-shaped intersections --marked with a cross in Figure \ref{fig:graph} (left)-- are given by the set of indices
$$\mathcal{I}_0^X=\{(2,3,4,8),(5,6,7,8),(1,6,9,10)\}.$$
Note that the set of indices $(4,5,6,9)$ also defines a closed path in the graph, but $(4,5,9), (5,6,9)\in \mathcal{I}_0^T$, so it does not represent an $X$-shaped intersection. Something similar applies to the sets of indices $(1,2,6,7)$, $(4,5,8,9)$ and $(1,2,7,8)$.

Finally, let us denote by $\gamma(i)$ the set of all adjacent fractures to subdomain $\Omega_i$, for $i\in \mathcal{I}_2$. In particular, $\gamma(1)=\{\gamma_{1,2},\,\gamma_{1,6},\,\gamma_{1,7},\,\gamma_{1,10}\}$, $\gamma(2)=\{\gamma_{1,2},\,\gamma_{2,3},\,\gamma_{2,7},\,\gamma_{2,8}\}$, and so on.
In turn, $\sigma(i,j)$ is defined as the set of all intersecting points in which the fracture $\gamma_{i,j}$ is involved, for $(i,j)\in \mathcal{I}_1$. In this case, $\sigma(1,2)=\{\sigma_{1,2,7}\}$, $\sigma(1,6)=\{\sigma_{1,6,7},\,\sigma_{1,6,9,10}\}$, and so on.

With the aim of defining an interface problem, each fracture permeability $\mathbf{K}_{i,j}$ is decomposed in a similar way to that introduced in \eqref{K:decomposition}, with corresponding coefficients $K_{i,j}^{\boldsymbol{\tau}}$ and $K_{i,j}^{\mathbf{n}}$. In this framework, the following mixed-dimensional problem is formulated, for $i\in \mathcal{I}_2$ and $(i,j)\in \mathcal{I}_1$ \cite{ami:ker:mar:rob:2006}, 
\begin{equation}\label{interface:problem:intersections}
\begin{alignedat}{2}
\mathbf{u}_i&=-\mathbf{K}_i\nabla p_i\hspace*{1.5cm}&&\hbox{in }\Omega_i, \\
\nabla\cdot\mathbf{u}_i&=q_i\hspace*{1.5cm}&&\hbox{in }\Omega_i,  \\
\mathbf{u}_{i,j}&=-d_{i,j} \,K_{i,j}^{\boldsymbol{\tau}}\, \nabla^{\boldsymbol{\tau}}_{i,j}\, p_{i,j}\hspace*{1.5cm}&&\hbox{in }\gamma_{i,j}, \\
\nabla^{\boldsymbol{\tau}}_{i,j}\cdot\mathbf{u}_{i,j}&=q_{i,j}+(\mathbf{u}_i\cdot\mathbf{n}_i+\mathbf{u}_j\cdot\mathbf{n}_j)\hspace*{1.5cm}&&\hbox{in }\gamma_{i,j}, \\
\alpha_{i,j}(p_i-p_{i,j})&=\xi\,\mathbf{u}_i\cdot\mathbf{n}_i-(1-\xi)\,\mathbf{u}_{j}\cdot\mathbf{n}_{j}\hspace*{1.5cm}&&\hbox{in }\gamma_{i,j}, \\
\alpha_{i,j}(p_j-p_{i,j})&=\xi\,\mathbf{u}_j\cdot\mathbf{n}_j-(1-\xi)\,\mathbf{u}_{i}\cdot\mathbf{n}_{i}\hspace*{1.5cm}&&\hbox{in }\gamma_{i,j}, \\
p_i&=0\hspace*{1.5cm}&&\hbox{on }\Gamma_i, \\
p_{i,j}&=0\hspace*{1.5cm}&&\hbox{on }\Gamma_{i,j},  
\end{alignedat}
\end{equation}
where $\alpha_{i,j}=2K_{i,j}^{\mathbf{n}}/d_{i,j}$ and $\Gamma_{i,j}=\partial\gamma_{i,j}\cap\Gamma$. The notations $\nabla^{\boldsymbol{\tau}}_{i,j}\cdot$ and $\nabla^{\boldsymbol{\tau}}_{i,j}$ stand for the tangential divergence and gradient operators, as defined by \eqref{tangential:div:grad}, on the fracture $\gamma_{i,j}$. At the intersections, we shall impose mass conservation and pressure continuity. In particular, at every $T$-shaped intersecting point $\sigma_{i,j,k}$, with $(i,j,k)\in \mathcal{I}_0^T$, we impose 
\begin{equation}\label{inters:conditions:t}
\begin{alignedat}{2}
\sum_{m,n\in\{i,j,k\},\,(m,n)\in\mathcal{I}_1}\mathbf{u}_{m,n}\cdot\mathbf{n}_{m,n}&=0,&&\\
p_{m,n}&=p_{i,j,k}\hspace*{0.7cm}&& \forall\,m,n\in\{i,j,k\},\,(m,n)\in\mathcal{I}_1. 
\end{alignedat}
\end{equation}
In turn, at every $X$-shaped intersecting point $\sigma_{i,j,k,l}$, with $(i,j,k,l)\in \mathcal{I}_0^X$, we impose 
\begin{equation}\label{inters:conditions:x}
\begin{alignedat}{2}
\sum_{m,n\in\{i,j,k,l\},\,(m,n)\in\mathcal{I}_1}\mathbf{u}_{m,n}\cdot\mathbf{n}_{m,n}&=0,&&\\
p_{m,n}&=p_{i,j,k,l}\hspace*{0.7cm}&& \forall\, m,n\in\{i,j,k,l\},\,(m,n)\in \mathcal{I}_1. 
\end{alignedat}
\end{equation} 
For a discussion on more general compatibility conditions at the intersections, we refer the reader to \cite{fle:fum:sco:16,for:fum:sco:ruf:14,fum:kei:18,sch:fle:hel:woh:15}.


In order to define the weak formulation of problem \eqref{interface:problem:intersections}-\eqref{inters:conditions:x}, we introduce the spaces for the velocity unknowns \cite{ami:ker:mar:rob:2006}
\begin{align*}
	&\mathbf{W}^2=\bigoplus_{i\in \mathcal{I}_2} \{\mathbf{v}_i\in H(\mathrm{div},\Omega_i): \mathbf{v}_i\cdot \mathbf{n}_i\in L^2(\gamma)\ \forall\,\gamma\in\gamma(i)\},\\[1ex]
	&\textbf{W}^1=\bigoplus_{(i,j)\in \mathcal{I}_1} \{\mathbf{v}_{i,j}\in H(\mathrm{div}_{i,j}^{\boldsymbol{\tau}},\gamma_{i,j}): \mathbf{v}_{i,j}\cdot \mathbf{n}_{i,j}\in L^2(\sigma)\  \forall\,\sigma\in\sigma(i,j)\},
\end{align*}
together with the spaces for the pressures
\begin{align*}
	&M^2=\bigoplus_{i\in \mathcal{I}_2}L^2(\Omega_i),&& M^1=\bigoplus_{(i,j)\in \mathcal{I}_1}L^2(\gamma_{i,j}),\\[1ex]
	&M^{0,T}=\bigoplus_{(i,j,k)\in \mathcal{I}_0^T}L^2(\sigma_{i,j,k}), && M^{0,X}=\bigoplus_{(i,j,k,l)\in \mathcal{I}_0^X}L^2(\sigma_{i,j,k,l}).
\end{align*}
Note that the superscript notation of these spaces provides information about the dimensionality of the corresponding domain in which they are defined. This dimensional decomposition framework has been proposed in \cite{boo:nor:yot:18}. Then, the global spaces $\textbf{W}=\textbf{W}^2\oplus\mathbf{W}^1$ and $M=M^2\oplus M^1\oplus M^{0,T}\oplus M^{0,X}$ are endowed with the norms
\begin{align*}
	\|\mathbf{v}\|^2_{\mathbf{W}}&=\sum_{i\in \mathcal{I}_2}\left(\|\mathbf{v}_i\|^2_{L^2(\Omega_i)}+\|\nabla\cdot\mathbf{v}_i\|^2_{L^2(\Omega_i)}+\sum_{\gamma\in\gamma(i)}\|\mathbf{v}_i\cdot\mathbf{n}_{i}\|^2_{L^2(\gamma)}\right)+\\[1ex]
	&\hspace*{0.7cm}\sum_{(i,j)\in \mathcal{I}_1}\left(\|\mathbf{v}_{i,j}\|^2_{L^2(\gamma_{i,j})}+\|\nabla^{\boldsymbol{\tau}}_{i,j}\cdot\mathbf{v}_{i,j}\|^2_{L^2(\gamma_{i,j})}+\sum_{\sigma\in\sigma(i,j)}\|\mathbf{v}_{i,j}\cdot\mathbf{n}_{i,j}\|^2_{L^2(\sigma)}\right),\\[3ex]
	\|r\|^2_M&=\sum_{i\in \mathcal{I}_2}\|r_i\|^2_{L^2(\Omega_i)}+\sum_{(i,j)\in \mathcal{I}_1}\|r_{i,j}\|^2_{{L^2(\gamma_{i,j})}}+\sum_{(i,j,k)\in \mathcal{I}_0^T}\|r_{i,j,k}\|^2_{L^2(\sigma_{i,j,k})}
	+\sum_{(i,j,k,l)\in \mathcal{I}_0^X}\|r_{i,j,k,l}\|^2_{L^2(\sigma_{i,j,k,l})}.
\end{align*}

In this framework, the bilinear forms $a:\mathbf{W}\times\mathbf{W}\rightarrow\mathbb{R}$ and $b:\mathbf{W}\times M\rightarrow\mathbb{R}$ are given, respectively, by
%
%
\begin{align*}
	&a(\mathbf{u},\mathbf{v})=\sum_{i\in \mathcal{I}_2}\left(\mathbf{K}_i^{-1}\mathbf{u}_i,\mathbf{v}_i\right)_{\Omega_i}+\sum_{(i,j)\in \mathcal{I}_1}\left((d_{i,j} \,K_{i,j}^{\boldsymbol{\tau}})^{-1}\mathbf{u}_{i,j},\mathbf{v}_{i,j}\right)_{\gamma_{i,j}}\nonumber\\
	&+\sum_{(i,j)\in \mathcal{I}_1}\left(\alpha_{i,j}^{-1}(\xi\,\mathbf{u}_i\cdot\mathbf{n}_i-(1-\xi)\,\mathbf{u}_{j}\cdot\mathbf{n}_{j}),\mathbf{v}_i\cdot\mathbf{n}_i\right)_{\gamma_{i,j}}\nonumber
	+\sum_{(i,j)\in \mathcal{I}_1}\left(\alpha_{i,j}^{-1}(\xi\,\mathbf{u}_{j}\cdot\mathbf{n}_{j}-(1-\xi)\,\mathbf{u}_{i}\cdot\mathbf{n}_{i}),\mathbf{v}_{j}\cdot\mathbf{n}_{j}\right)_{\gamma_{i,j}},
\end{align*}
and
\begin{align*}
	&b(\mathbf{u},r)=\sum_{i\in \mathcal{I}_2}\left(\nabla\cdot \mathbf{u}_i,r_i\right)_{\Omega_i}+\sum_{(i,j)\in \mathcal{I}_1}\left(\nabla^{\boldsymbol{\tau}}_{i,j}\cdot \mathbf{u}_{i,j},r_{i,j}\right)_{\gamma_{i,j}}
	-\sum_{(i,j)\in \mathcal{I}_1}\left(\mathbf{u}_i\cdot\mathbf{n}_i+\mathbf{u}_{j}\cdot\mathbf{n}_{j},r_{i,j}\right)_{\gamma_{i,j}}\nonumber\\
	&-\!\!\!\sum_{(i,j,k)\in \mathcal{I}_0^T}\!\!\!\left(\sum_{(m,n)\in\{i,j,k\},\,(m,n)\in \mathcal{I}_1}\!\!\!\mathbf{u}_{m,n}\cdot\mathbf{n}_{m,n},r_{i,j,k}\!\!\right)_{\!\!\!\sigma_{i,j,k}}
	-\!\!\!\sum_{(i,j,k,l)\in \mathcal{I}_0^X}\!\!\!\left(\sum_{(m,n)\in\{i,j,k,l\},\,(m,n)\in \mathcal{I}_1}\!\!\!\mathbf{u}_{m,n}\cdot\mathbf{n}_{m,n},r_{i,j,k,l}\!\!\right)_{\!\!\!\sigma_{i,j,k,l}}.
\end{align*}
In turn, the linear form $L:M\rightarrow\mathbb{R}$ associated with the source terms is defined as
$$
L(r)=\sum_{i\in \mathcal{I}_2}\left(q_i,r_i\right)_{\Omega_i}+\sum_{(i,j)\in \mathcal{I}_1}\left(q_{i,j},r_{i,j}\right)_{\gamma_{i,j}}.
$$
%
In this setting, the weak formulation of problem \eqref{interface:problem:intersections}-\eqref{inters:conditions:x} shows the same structure as \eqref{weak:form}. However, since the function spaces and forms are newly defined in this case, we reproduce it here for convenience: \emph{Find} $(\mathbf{u},p)\in\mathbf{W}\times M$ \emph{such that}
\begin{equation}\label{weak:form:inters}
\begin{alignedat}{2}
a(\mathbf{u},\mathbf{v})-b(\mathbf{v},p)&=0 \hspace*{1.5cm}&& \forall\, \mathbf{v}\in\mathbf{W},\\
b(\mathbf{u},r)&=L(r)\hspace*{1.5cm} && \forall r\in M.
\end{alignedat}
\end{equation}
Following \cite{ami:ker:mar:rob:2006,Martin-Jaffre-Roberts}, if the permeabilities in the subdomains and the coefficients $K_{i,j}^{\mathbf{n}}/d_{i,j}$ and $K_{i,j}^{\boldsymbol{\tau}}d_{i,j}$ are bounded by positive constants, the existence and uniqueness of solution of this problem can be proved for the case $\xi>1/2$.

In the remaining of this subsection, we introduce a mixed finite element discretization for problem \eqref{weak:form:inters}. Along the lines of Subsection \ref{subsec:mfe}, we suppose that the subdomains $\Omega_i$ admit a rectangular partition $\mathcal{T}^h_{i}$, for $i\in\mathcal{I}_2$. Further, the meshes $\mathcal{T}^h_{i}$ are assumed to match at the interfaces $\gamma_{i,j}$, i.e., they induce a unique partition $\mathcal{T}^h_{i,j}$ on $\gamma_{i,j}$, for $(i,j)\in \mathcal{I}_1$. Note that the intersecting points of the fractures $\{\sigma_{i,j,k}\}_{(i,j,k)\in\mathcal{I}_0^T}$ and $\{\sigma_{i,j,k,l}\}_{(i,j,k,l)\in\mathcal{I}_0^X}$ are vertices of some of the preceding meshes. For instance, given $(i,j,k)\in\mathcal{I}_0^T$, the intersecting point $\sigma_{i,j,k}$ is a vertex of the two-dimensional meshes $\mathcal{T}^h_{i}$, $\mathcal{T}^h_{j}$ and $\mathcal{T}^h_{k}$, and also a vertex of the one-dimensional meshes $\mathcal{T}^h_{i,j}$, $\mathcal{T}^h_{j,k}$ and $\mathcal{T}^h_{i,k}$.


In this context, let $\mathbf{W}^h_{i}\times M^h_{i}$ be the lowest order Raviart--Thomas mixed finite element spaces defined on $\mathcal{T}^h_{i}$, for $i\in\mathcal{I}_2$, let $\mathbf{W}^h_{{i,j}}\times M^h_{{i,j}}$ be the lowest order Raviart--Thomas mixed finite element spaces defined on $\mathcal{T}^h_{{i,j}}$, for $(i,j)\in\mathcal{I}_1$, and let $M^h_{{i,j,k}}$ and $M^h_{{i,j,k,l}}$ be equal to $\mathbb{R}$, for $(i,j,k)\in\mathcal{I}_0^T$ and $(i,j,k,l)\in\mathcal{I}_0^X$. Based on the notations
\begin{align*}
&\mathbf{W}^{h,2}=\displaystyle\bigoplus_{i\in \mathcal{I}_2}\mathbf{W}^h_{i},&&\mathbf{W}^{h,1}=\displaystyle\bigoplus_{(i,j)\in \mathcal{I}_1}\mathbf{W}^h_{{i,j}},\\[1ex]
&M^{h,2}=\displaystyle\bigoplus_{i\in \mathcal{I}_2}M^h_{i},&& M^{h,1}=\displaystyle\bigoplus_{(i,j)\in \mathcal{I}_1}M^h_{{i,j}},\\[1ex]
& M^{h,0,T}=\displaystyle\bigoplus_{(i,j,k)\in \mathcal{I}_0^T}M^h_{{i,j,k}},&& M^{h,0,X}=\displaystyle\bigoplus_{(i,j,k,l)\in \mathcal{I}_0^X}M^h_{{i,j,k,l}},
\end{align*}
the approximation spaces are obtained as $\mathbf{W}^h=\mathbf{W}^{h,2}\oplus\mathbf{W}^{h,1}$ for the velocities, and $M^h=M^{h,2}\oplus M^{h,1}\oplus M^{h,0,T} \oplus M^{h,0,X}$ for the pressures. Then, using the quadrature rules introduced in the previous subsection, we define the discrete forms
\begin{align*}
&a_{h}(\mathbf{u},\mathbf{v})=\sum_{i\in \mathcal{I}_2}\left(\mathbf{K}_i^{-1}\mathbf{u}_i,\mathbf{v}_i\right)_{\Omega_i,\mathbf{TM}}+\sum_{(i,j)\in \mathcal{I}_1}\left((d_{i,j} \,K_{i,j}^{\boldsymbol{\tau}})^{-1}\mathbf{u}_{i,j},\mathbf{v}_{i,j}\right)_{\gamma_{i,j},\mathbf{T}}\nonumber\\
	&+\!\!\sum_{(i,j)\in \mathcal{I}_1}\!\!\left(\alpha_{i,j}^{-1}(\xi\,\mathbf{u}_i\cdot\mathbf{n}_i-(1-\xi)\,\mathbf{u}_{j}\cdot\mathbf{n}_{j}),\mathbf{v}_i\cdot\mathbf{n}_i\right)_{\gamma_{i,j}}
	+\!\!\sum_{(i,j)\in \mathcal{I}_1}\!\!\left(\alpha_{i,j}^{-1}(\xi\,\mathbf{u}_{j}\cdot\mathbf{n}_{j}-(1-\xi)\,\mathbf{u}_{i}\cdot\mathbf{n}_{i}),\mathbf{v}_{j}\cdot\mathbf{n}_{j}\right)_{\gamma_{i,j}},\nonumber\\[1ex]
	&L_h(r)=\sum_{i\in \mathcal{I}_2}\left(q_i,r_i\right)_{\Omega_i,\mathbf{M}}+\sum_{(i,j)\in \mathcal{I}_1}\left(q_{i,j},r_{i,j}\right)_{\gamma_{i,j},\mathbf{M}}.
\end{align*}
Finally, the mixed finite element approximation to \eqref{weak:form:inters} takes the form \eqref{mfe:approx} for the newly defined discrete spaces and forms, namely: \emph{Find} $(\mathbf{u}_h,p_h)\in \mathbf{W}^h\times M^h$ \emph{such that}
%
\begin{equation}\label{mfe:approx:inters}
\begin{alignedat}{2}
a_{h}(\mathbf{u}_h,\mathbf{v}_h)-b(\mathbf{v}_h,p_h)&=0&&\hspace*{1cm}\forall\,\mathbf{v}_h\in\mathbf{W}^h,\\
b(\mathbf{u}_h,r_h)&=L_h(r_h)&&\hspace*{1cm}\forall\,r_h\in M^h.
\end{alignedat}
\end{equation}

In order to derive the algebraic linear system underlying \eqref{mfe:approx:inters}, we group the pressure and velocity unknowns taking into account their dimensions, i.e., $\mathbf{u}_h=[U^2,U^1]^T$ and $p_h=[P^2,P^1,P^0]^T$. Note that $U^2$ and $U^1$ are vectors related to the velocity unknowns on the two-dimensional subdomains and the one-dimensional fractures, respectively.
In turn, $P^2$, $P^1$ and $P^0$ are vectors related to the pressure unknowns on the two-dimensional subdomains, the one-dimensional fractures and the zero-dimensional intersections, respectively. In this context, the algebraic linear system stemming from \eqref{mfe:approx:inters} may be written as the saddle point problem
\begin{equation}\label{whole_system}
\begin{bmatrix}
A_{2,2} & 0 & B_{2,2}^T & F_{2,1}^T & 0\\[0.5ex]
0 & A_{1,1} & 0 & B_{1,1}^T & F_{1,0}^T\\[0.5ex]
B_{2,2} & 0 & 0 & 0 & 0 \\[0.5ex]
F_{2,1} & B_{1,1} & 0 & 0 & 0 \\[0.5ex]
0 & F_{1,0} & 0 & 0 & 0
\end{bmatrix}
\begin{bmatrix}
U^2\\[0.5ex]U^1\\[0.5ex]P^2\\[0.5ex]P^1\\[0.5ex]P^0
\end{bmatrix}=
\begin{bmatrix}
0\\[0.5ex]0\\[0.5ex]Q^2\\[0.5ex]Q^1\\[0.5ex]0
\end{bmatrix}.
\end{equation}
This is a generalization of the linear system \eqref{linear:system:1fracture} obtained in the case of a single fracture. In such a case, $U^2$ and $P^2$ were composed of two blocks, one per subdomain, while $U^1$ and $P^1$ consisted of one block corresponding to the only fracture $\gamma$, and $P^0$ was lacking since there were no intersecting points. In contrast, for the example shown in Figure \ref{fig:graph}, $U^2$ and $P^2$ would be composed of 10 blocks, one per subdomain, $U^1$ and $P^1$ would consist of 18 blocks, one per fracture, and $P^0$ would group 9 pressure unknowns, one per intersecting point of fractures.

\section{Monolithic multigrid methods for mixed-dimensional elliptic problems}\label{sec:3}
\setcounter{section}{3}


In this work, we propose an efficient monolithic multigrid solver for flow in fractured porous media. Multigrid methods \cite{Bra77,TOS01} aim to accelerate the slow convergence of classical iterative methods by using coarse meshes. Since these latter have a strong smoothing effect on the error of the solution, this error can be properly represented in coarser grids where the computations are much less expensive.

Suppose that $A_k \, u_k = f_k$ is the system to solve, where the matrix $A_k$ corresponds to a discretization of a partial differential equation on a grid $G^k$, $f_k$ is the right-hand side and $u_k$ is the unknown vector. In order to apply a standard two-grid cycle for solving this problem, we perform the following steps:\\ 

\begin{enumerate}
\item Apply $\nu_1$ iterations of a classical iterative method, called smoother, on $G^k$ (pre-smoothing step). 
\item Compute the residual of the current fine grid approximation.
\item Restrict the residual to the coarse grid $G^{k-1}$ by using a restriction operator $R_{k}^{k-1}$.
\item Solve the residual equation on the coarse grid.
\item Interpolate the obtained correction to the fine grid $G^k$ by using a prolongation operator $P_{k-1}^{k}$.
\item Add the interpolated correction to the current  fine grid approximation. 
\item Apply $\nu_2$ iterations of a classical iterative method on $G^k$ (post-smoothing step).\\
\end{enumerate}

Since we do not need to solve the problem on the coarse grid exactly, we can apply the same algorithm in a recursive way by using a hierarchy of coarser meshes, giving rise to the well-known multigrid cycle. It is clear that many details have to be fixed for the design of an efficient multigrid method, since all the components have to be properly chosen. In particular, we need to specify the hierarchy of grids, the coarse-grid operators, the type of cycle, the inter-grid transfer operators and the smoothing procedure. Next, we explain our choices in this work. 

\subsection{Hierarchy of meshes, coarse-grid operators and cycle type}

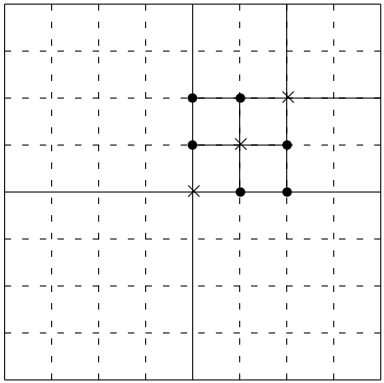
\begin{figure}[t]
	\begin{center}
		\begin{minipage}{0.4\textwidth}
			\begin{center}
				\unitlength=0.05cm
				\begin{picture}(100,100)
				\put(0,0){\line(1,0){100}}\put(0,0){\line(0,1){100}}\put(100,0){\line(0,1){100}}
				\put(0,100){\line(1,0){100}}
				\dashline{2}(12.5,0)(12.5,100)
				\dashline{2}(25,0)(25,100)
				\dashline{2}(37.5,0)(37.5,100)
				\put(50,0){\line(0,1){100}}
				\dashline{2}(62.5,0)(62.5,100)
				\dashline{2}(75,0)(75,100)
				\dashline{2}(87.5,0)(87.5,100)
				\put(0,50){\line(1,0){100}}
				\put(50,75){\line(1,0){50}}
				\put(75,50){\line(0,1){50}}
				\put(50,62.5){\line(1,0){25}}
				\put(62.5,50){\line(0,1){25}}
				\dashline{2}(0,12.5)(100,12.5)
				\dashline{2}(0,25)(100,25)
				\dashline{2}(0,37.5)(100,37.5)
				\dashline{2}(0,62.5)(100,62.5)
				\dashline{2}(0,75)(100,75)
				\dashline{2}(0,87.5)(100,87.5)
				\put(48.4,61.05){\footnotesize$\bullet$}
				\put(48.4,73.55){\footnotesize$\bullet$}
				\put(61.12,73.55){\footnotesize$\bullet$}
				\put(61.1,48.5){\footnotesize$\bullet$}
				\put(73.5,48.5){\footnotesize$\bullet$}
				\put(73.5,61.05){\footnotesize$\bullet$}
				\put(47.8,48.6){\footnotesize$\times$}
				\put(60.3,61.1){\footnotesize$\times$}
				\put(72.9,73.6){\footnotesize$\times$}
				\end{picture}
			\end{center}
		\end{minipage}
	\end{center}
	\caption{Coarsest grid corresponding to the fracture network shown in Figure \ref{fig:graph}.}
	\label{fig:graph_grid}
\end{figure}

The implementation of a geometric multigrid method requires to define the problem on grids of various sizes, namely a hierarchy of grids. Here, such a hierarchy is constructed in the following way. First, we consider a coarse grid which is built taking into account the location of the fractures.  This mesh is generated by assuming that every fracture coincides with an edge of some element in the porous medium grid. As explained in Section \ref{sec:2}, we suppose that the grids in the subdomains match at the interfaces. Thus, in the case of considering $m$ subdomains, we define $G^0 = \textstyle\bigcup_{i=1}^m {\mathcal T}_i^h$ as the coarsest possible grid fulfilling this criterion. As an example, if we consider the fracture configuration shown in Figure \ref{fig:graph}, the coarsest grid $G^0$ is given in Figure \ref{fig:graph_grid}.

\begin{figure}[t]
	\begin{center}
		\begin{minipage}{0.25\textwidth}
			\begin{center}
				\unitlength=0.02cm
				\begin{picture}(150,300) 
					\put(0,25){\line(1,0){150}}
					\thicklines
					\put(0,125){\line(1,0){150}}
					\thinlines
					\put(0,225){\line(1,0){150}}
					\put(25,0){\line(0,1){250}} 
					\put(125,0){\line(0,1){250}}	
					\put(69.2,72.2){\footnotesize$\times$}	
					\put(21,75){\line(1,0){8}}	
					\put(121,75){\line(1,0){8}}
					\put(75,21){\line(0,1){8}}	
					\put(75,121){\line(0,1){8}}	
					\put(69.2,172.2){\footnotesize$\times$}
					\put(21,175){\line(1,0){8}}	
					\put(121,175){\line(1,0){8}}
					\put(75,221){\line(0,1){8}}	
					\put(75,125){\circle{10}}
					\put(25,125){\circle*{8}}
					\put(125,125){\circle*{8}}
					\put(70,265){\footnotesize$G^{k}$} 
				\end{picture}
			\end{center}
		\end{minipage}	
		\begin{minipage}{0.05\textwidth}
			\begin{center}
				\vspace{1.1cm}
				{\LARGE $\Rightarrow$}
			\end{center}		
		\end{minipage}
		\begin{minipage}{0.25\textwidth}
			\begin{center}
				\unitlength=0.02cm
				\begin{picture}(150,300)
					\put(70,265){\footnotesize$G^{k+1}$}
					\put(0,25){\line(1,0){150}}
					\put(0,75){\line(1,0){150}}
					\thicklines
					\put(0,125){\line(1,0){150}}
					\thinlines
					\put(0,225){\line(1,0){150}}
					\put(0,175){\line(1,0){150}}
					\put(25,0){\line(0,1){250}}
					\put(75,0){\line(0,1){250}}
					\put(125,0){\line(0,1){250}}	
					\put(21,50){\line(1,0){8}}
					\put(71,50){\line(1,0){8}}	
					\put(121,50){\line(1,0){8}}				
					\put(21,100){\line(1,0){8}}	
					\put(71,100){\line(1,0){8}}
					\put(121,100){\line(1,0){8}}
					\put(21,150){\line(1,0){8}}
					\put(71,150){\line(1,0){8}}	
					\put(121,150){\line(1,0){8}}
					\put(21,200){\line(1,0){8}}
					\put(71,200){\line(1,0){8}}	
					\put(121,200){\line(1,0){8}}							
					\put(50,21){\line(0,1){8}}	
					\put(50,71){\line(0,1){8}}
					\put(50,121){\line(0,1){8}}
					\put(50,171){\line(0,1){8}}
					\put(50,221){\line(0,1){8}}															\put(100,21){\line(0,1){8}}	
					\put(100,71){\line(0,1){8}}
					\put(100,121){\line(0,1){8}}
					\put(100,171){\line(0,1){8}}
					\put(100,221){\line(0,1){8}}	
					\put(44.2,47.2){\footnotesize$\times$}	
					\put(44.2,97.2){\footnotesize$\times$}
					\put(44.2,147.2){\footnotesize$\times$}
					\put(44.2,197.2){\footnotesize$\times$}
					\put(94.2,47.2){\footnotesize$\times$}	
					\put(94.2,97.2){\footnotesize$\times$}
					\put(94.2,147.2){\footnotesize$\times$}	
					\put(94.2,197.2){\footnotesize$\times$}
					\put(75,125){\circle*{8}}
					\put(25,125){\circle*{8}}
					\put(125,125){\circle*{8}}
					\put(100,125){\circle{10}}
					\put(50,125){\circle{10}}
				\end{picture}
			\end{center}
		\end{minipage}\hspace*{-0.5cm}
		\begin{minipage}{0.4\textwidth}
			\vspace{1cm}
			\begin{itemize}
			\item[$\times$] {\footnotesize pressure in the porous matrix}
			\item[$+$] {\footnotesize velocity in the porous matrix}
			\item[$\circ$]{\footnotesize pressure in the fracture}
			\item[$\bullet$]{\footnotesize velocity in the fracture}
			\end{itemize}
		\end{minipage}
	\end{center}
	\caption{Grid refinement procedure and location of the unknowns for both porous matrix and fracture (in bold line).}
	\label{fig:refinement}
\end{figure}
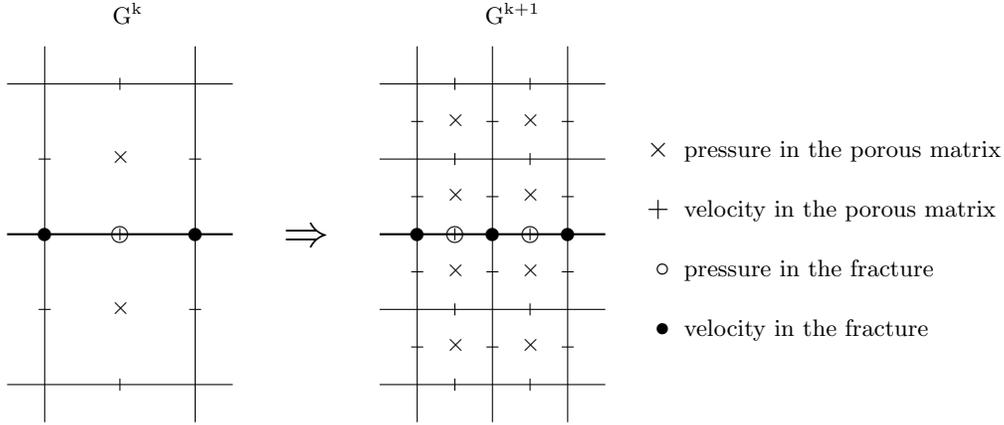

Then, the hierarchy of computational grids is created by applying a regular refinement process to each cell in that initial mesh. This means that we obtain a sequence of successively finer grids $G^0, G^1,\ldots, G^M$. In particular, since we are considering a quadrilateral partition of the porous medium, $G^{k+1}$ is obtained from $G^k$ by dividing each cell into four new elements for the next finer grid, as shown in Figure \ref{fig:refinement}, and this process continues until a fine enough target grid $G^M$ is obtained. 

Once the mesh hierarchy is generated, we consider a direct discretization of our problem on each grid. 
As for the type of cycle, we use W-cycles, since we have seen that this choice gives very good results for solving difficult coupled problems like the Darcy--Stokes system \cite{PeiyaoSISC} and the Biot--Stokes system \cite{PeiyaoJCP}.

\subsection{Inter-grid transfer operators}
Now, we define the restriction and interpolation operators involved in the multigrid method for solving the mixed-dimensional problem. We consider different transfer operators for the unknowns belonging to the matrix and for those located at the fractures. In particular, we choose two-dimensional and one-dimensional transfer operators, respectively. This means that we implement mixed-dimensional transfer operators in our multigrid algorithm in order to handle the problem at once. In matrix form, the chosen restriction operator $R_{k}^{k-1}$ from grid $G^k$ to $G^{k-1}$ is a block diagonal matrix since it does neither mix velocities and pressures nor unknowns in the porous matrix and in the fractures.

\begin{figure}[t]
\begin{tabular}{ccc}
\hspace{-0.5cm}\includegraphics[width = 0.35\textwidth]{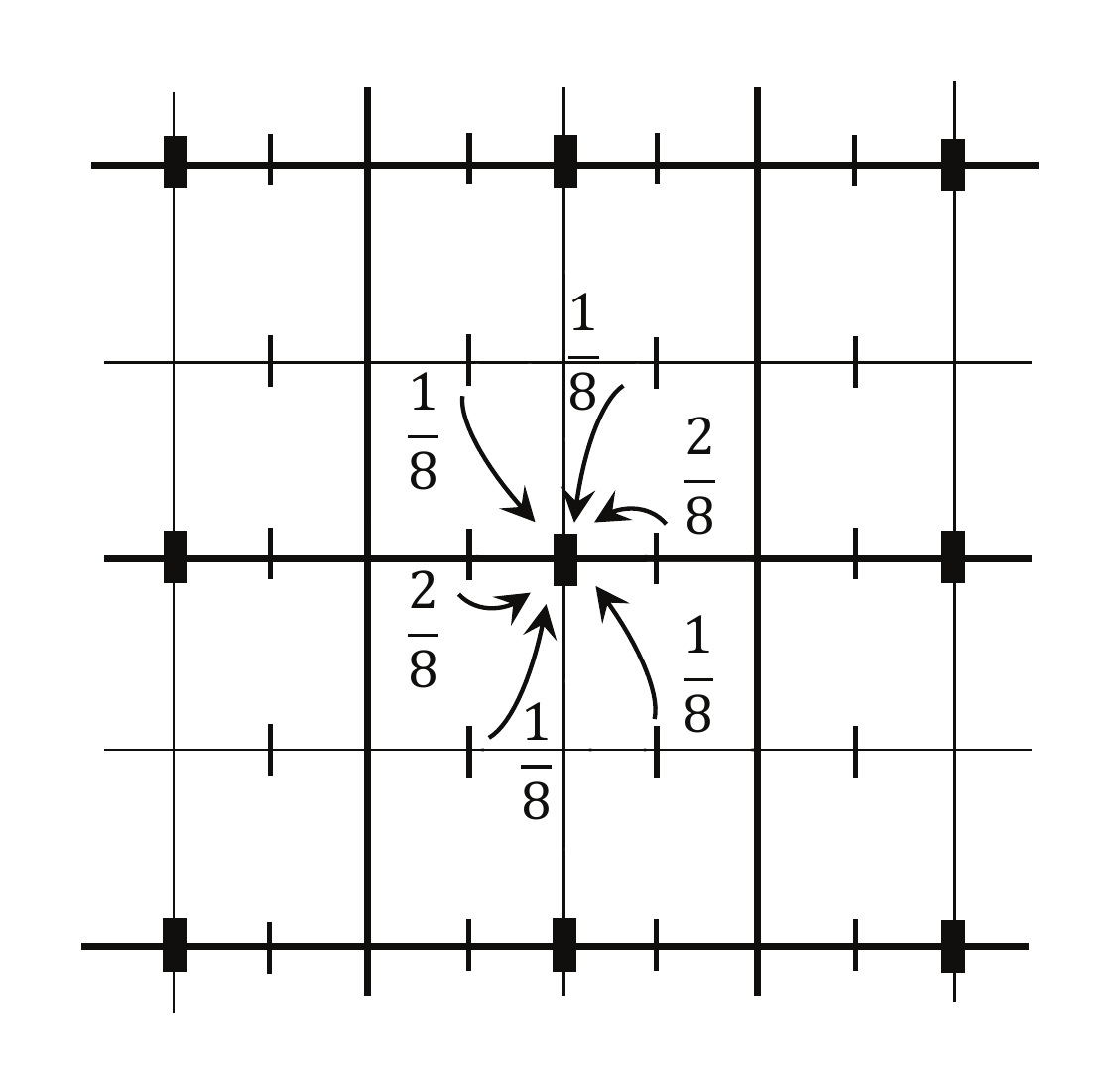}
&
\hspace{-0.6cm}\includegraphics[width = 0.35\textwidth]{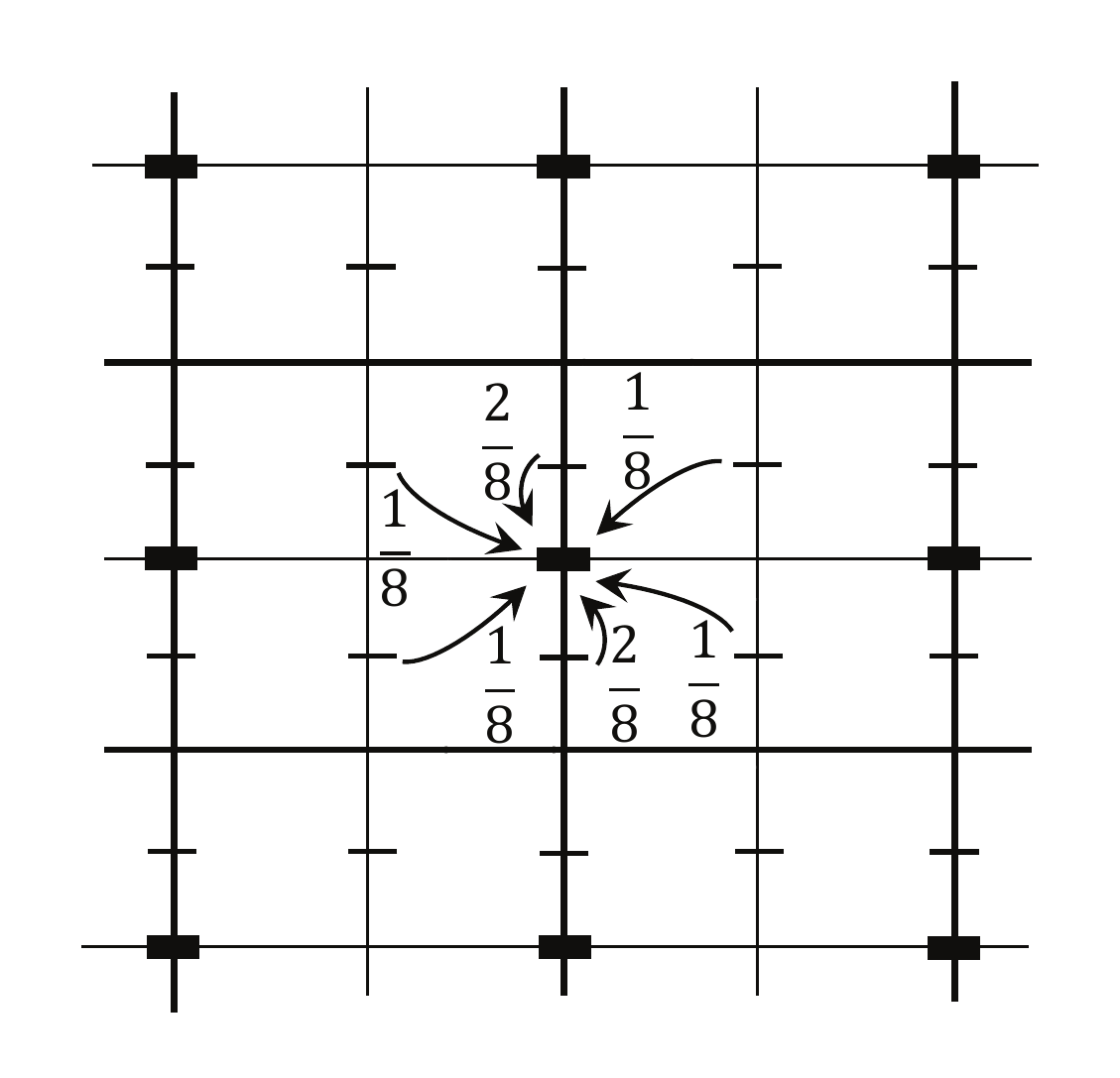}
&
\hspace{-0.6cm}\includegraphics[width = 0.35\textwidth]{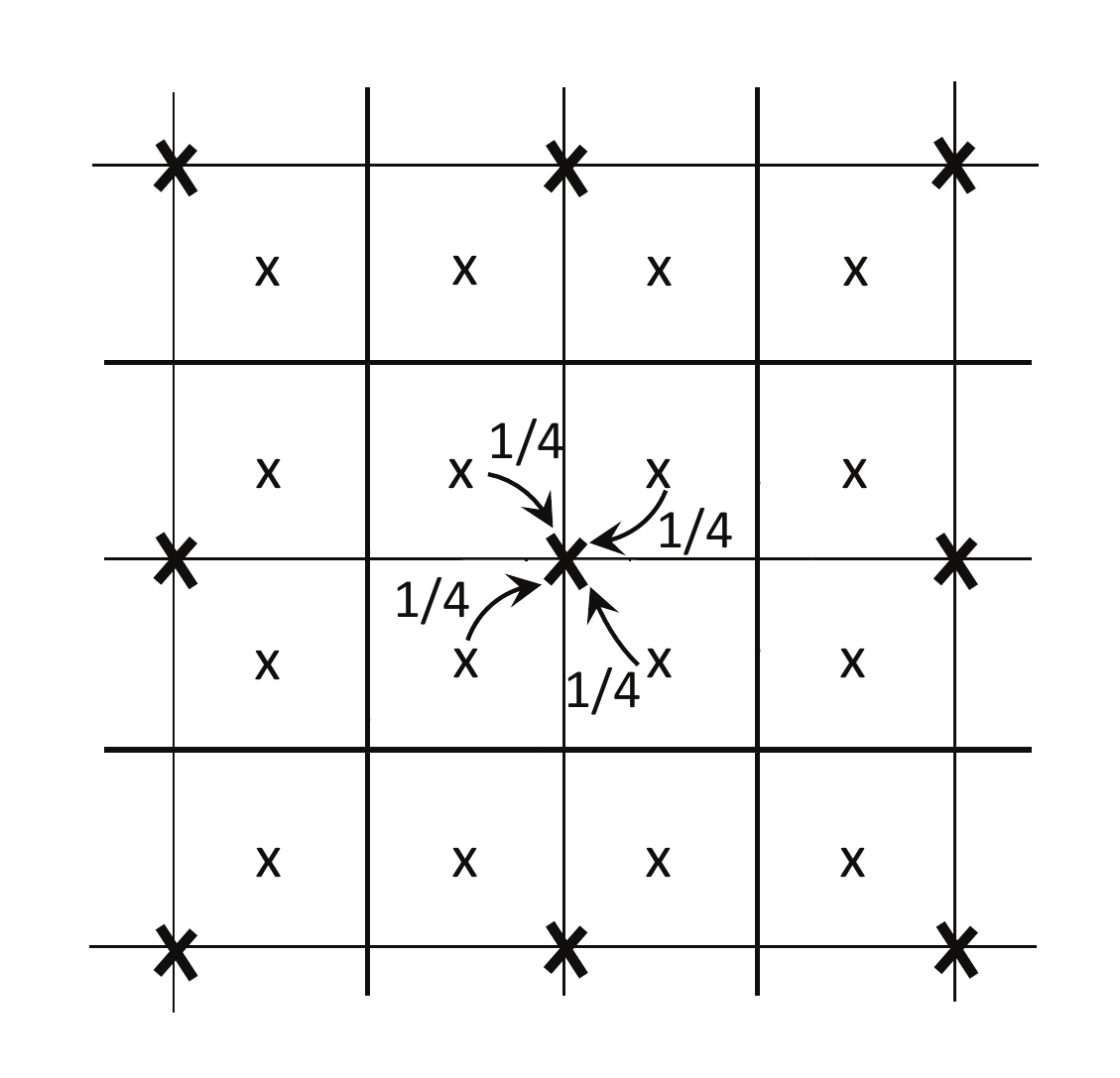}
\hspace{-0.6cm}
\end{tabular}
\label{fig:2Drestriction}
\caption{Restriction operators acting at the porous media unknowns}
\end{figure}

Due to the use of quadrature rules in this work, the mixed finite element method turns into a finite difference scheme on a staggered grid. As a consequence, we consider the standard restriction operators used for this type of meshes. Regarding the unknowns of the porous medium, we take into account the staggered arrangement of their location. Thus, the inter-grid transfer operators that act in the 
porous media unknowns are defined as follows: a six-point restriction is considered at velocity grid points, and a four-point restriction is applied at pressure grid points, as can be seen in Figure \ref{fig:2Drestriction}. 
The prolongation operator $P_{k-1}^{k}$, is chosen to be the adjoint of the restriction.

\begin{figure}[t]
\begin{center}
\begin{tabular}{cc}
\includegraphics[height = 2cm]{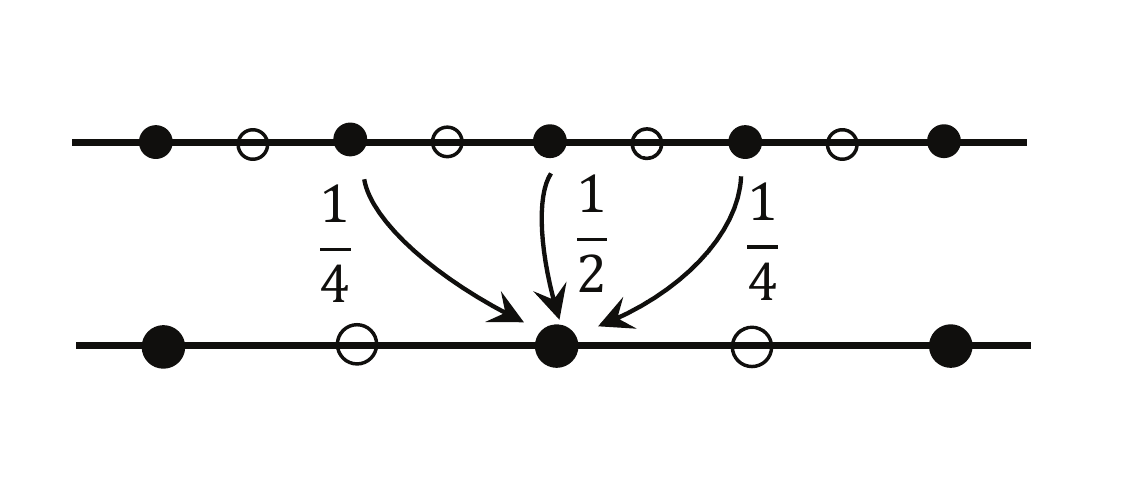}
&
\includegraphics[height = 2cm]{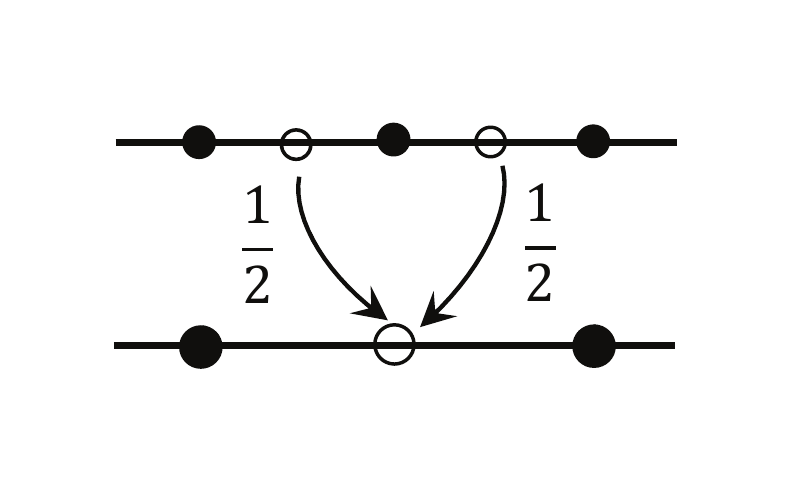}
\end{tabular}
\end{center}
\label{fig:1Drestriction}
\caption{Restriction operators acting at the fracture unknowns}
\end{figure}

Regarding the inter-grid transfer operators for the unknowns at the fractures, we again take into account 
their one-dimensional staggered arrangement, yielding the restriction transfer operators shown in Figure \ref{fig:1Drestriction}. 
Finally, the prolongation operators are chosen to be the corresponding adjoints.

\subsection{Smoother}
The performance of a multigrid method is essentially influenced by the smoothing algorithm. Here, in orderto deal with the difficulties generated by a saddle point problem, we consider a relaxation iteration among the class of multiplicative Schwarz smoothers. Basically, this type of iterations can be described as an overlapping block Gauss-Seidel method, where a small linear system of equations for each grid point has to be solved at each smoothing step. A particular case of such relaxation schemes is the so-called Vanka smoother, introduced in \cite{vanka} for solving the staggered finite difference discretization of the Navier--Stokes equations.

Due to the mixed-dimensional character of our problem, we propose a smoother $S_h$ which is written as the composition of three relaxation procedures acting on the two-dimensional cells of the porous matrix, $S_h^2$, the one-dimensional elements in the fractures, $S_h^1$, and the zero-dimensional intersection points, $S_h^0$, i.e.: $S_h = S_h^0 S_h^1 S_h^2$. Next, we describe these partial relaxation procedures:\\

\begin{enumerate}
\item\emph{Relaxation for the porous matrix.} 
The smoother considered for the unknowns located outside the fractured part of the domain is based on simultaneously updating all the unknowns appearing in the discrete divergence operator in the pressure equation. This way of building the blocks is very common in the Vanka-type smoothers used for Stokes and Navier-Stokes problems. This approach implies that four unknowns corresponding to velocities and one pressure unknown, see Figure \ref{vanka_blocks} (a), are relaxed simultaneously, making necessary to solve a $5 \times 5$ system for each cell. Then, we iterate over all the elements in lexicographic order, and for each of them the corresponding box is solved. \\


\begin{figure}
	\begin{center}
		\begin{tabular}{ccc}
			\centering\includegraphics[scale = 0.25]{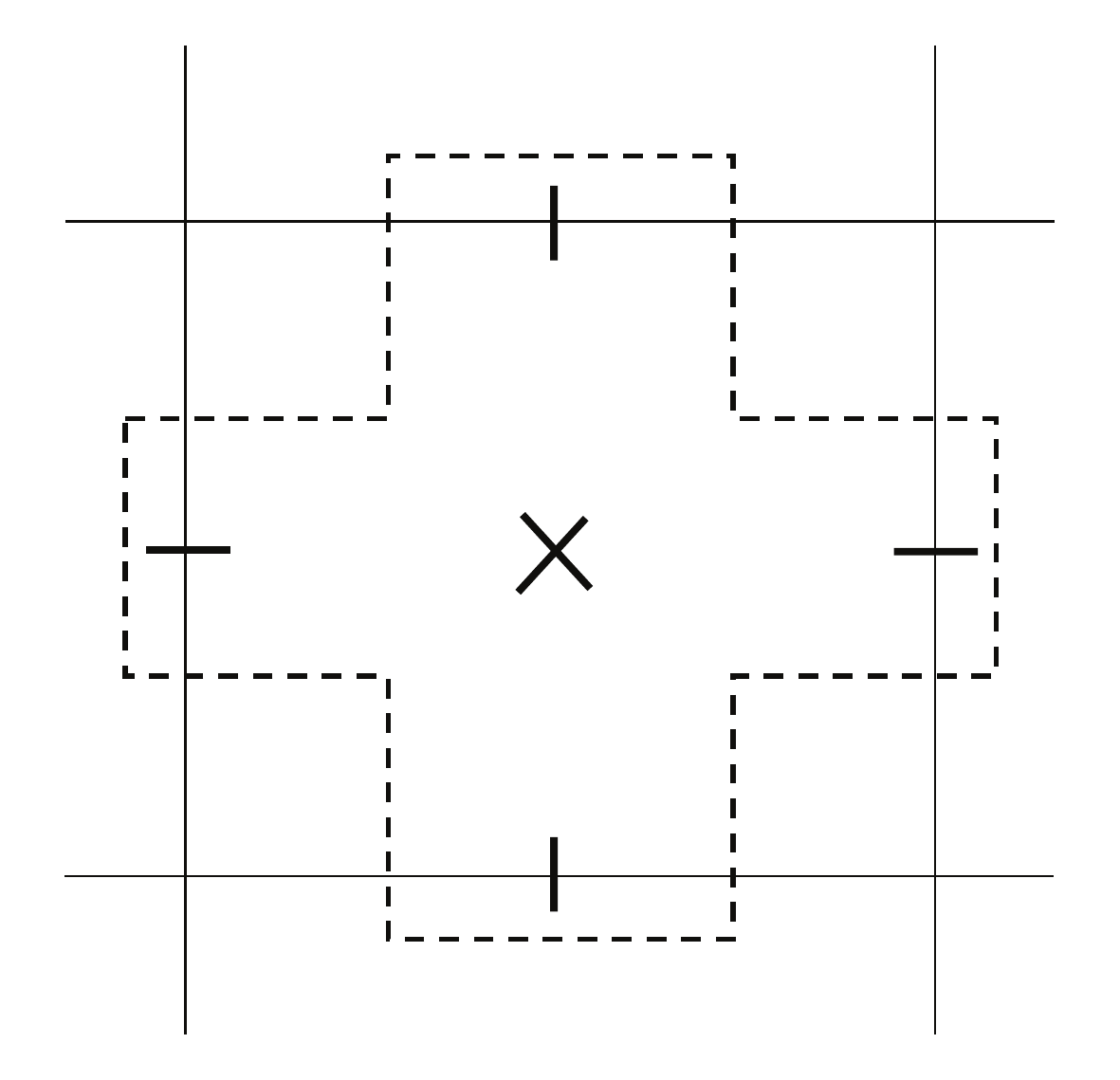}
			&
			\includegraphics[scale = 0.25]{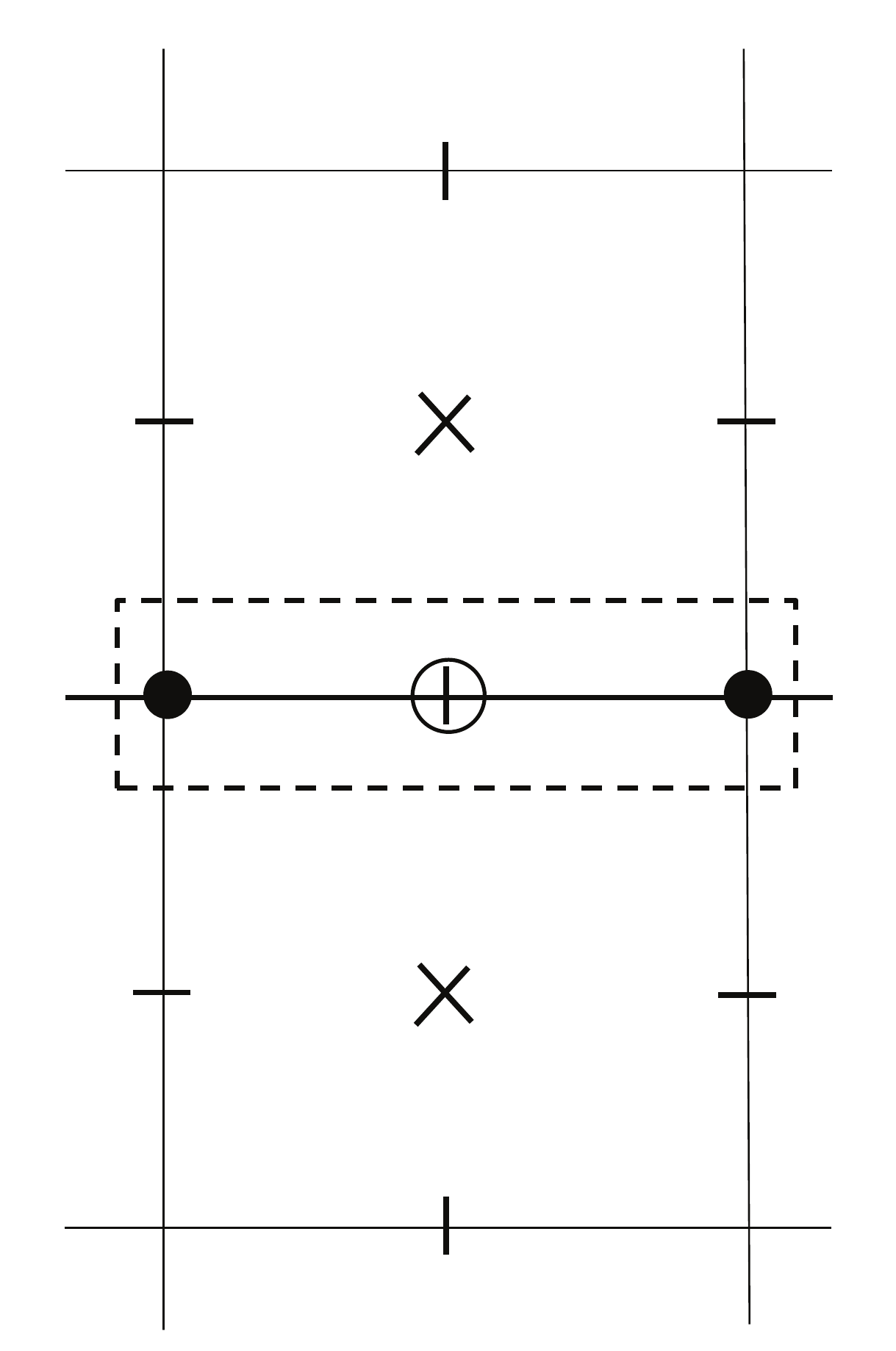}
			&
			\begin{minipage}{0.4\textwidth}
				\vspace{-4cm}
				\begin{itemize}
				\item[$\times$] {\footnotesize pressure in the porous matrix}
				\item[$+$] {\footnotesize velocity in the porous matrix}
				\item[$\circ$]{\footnotesize pressure in the fracture}
				\item[$\bullet$]{\footnotesize velocity in the fracture}
				\end{itemize}
			\end{minipage}\\
			(a) & (b) & 
		\end{tabular}
	\end{center}
	\caption{Unknowns updated together by the vanka-type smoothers applied (a) outside the fractures and (b) within the fractures.}
	\label{vanka_blocks}
\end{figure}

\item\emph{Relaxation for the fractures.}
The relaxation step applied to the unknowns located at the fractures is again based on simultaneously updating all the unknowns appearing in the discrete divergence operator in the pressure equation. This means that, in this case, for each element in the fracture we update five unknowns, three of them corresponding to the fracture and two of them to the matrix. In particular, each pressure unknown in the fracture is updated together with the two fracture velocities within the same element and the two porous matrix velocities located at the edges of the corresponding two-dimensional elements that match with that particular fracture element. This can be seen in Figure \ref{vanka_blocks} (b).
Notice that there are three unknowns located at the same point, the pressure in the fracture and
the two velocities corresponding to the elements adjacent to the fracture.\\

\item\emph{Relaxation for the intersections.}
At the intersection points of the fractures we apply a block Gauss--Seidel smoother coupling the fracture velocity unknowns located at each intersection, so that we need to solve a $2\times 2$, $3\times 3$ or $4\times 4$ system of equations on each of these grid points.\\

\end{enumerate}

The previously defined partial relaxation procedures can be formally written as
$$S_h^n= \prod_{B =1}^{N_{B,n}}\left(I-V_{B,n}^T (A^{B,n})^{-1} V_{B,n} A \right),\qquad\hbox{for } n=0,1,2,$$
where $A$ is the system matrix in \eqref{whole_system}, $N_{B,n}$ is the number of $n$-dimensional elements in the partition, $V_{B,n}$ represents the projection operator from the unknown vector to the vector of unknowns involved in the block to solve, and matrix $A^{B,n}$ is defined as $A^{B,n} = V_{B,n} A V_{B,n}^T$.

\subsection{Implementation}

The proposed monolithic mixed-dimensional multigrid method is implemented in a blockwise manner. Given an arbitrary fracture network composed of vertical and horizontal fractures, the first step is to construct  a uniform rectangular coarse grid so that the fracture network is contained in the set of edges of the grid. After that, a regular refinement process is applied on each block in the coarse grid until a target mesh with an appropriate fine grid scale to solve the problem is obtained. Then, for each step of the multigrid method, the two-dimensional components are performed in the porous matrix grid points whereas within the fractures one-dimensional components are implemented (notice that in the smoother, for example, this one-dimensional computation includes also unknowns from the porous matrix). In particular, in the smoother first the unknowns in the porous matrix are relaxed by using the standard two-dimensional Vanka smoother for Darcy problem, and after that, a one-dimensional Vanka smoother is used to update the unknowns located within the fractures. Finally, at the intersection points between different fractures, the velocities from different fractures are simultaneously relaxed.  
\begin{remark}
	Notice that this strategy can be easily extended to triangular grids in order to deal with more complex fracture networks. The idea would be to construct a non-structured coarse triangulation in such a way that the fracture network is part of its edges, and then to apply a regular refinement to the input triangles in order to obtain a semi-structured triangular grid in which the geometric multigrid method can be easily applied (see \cite{myBook, ARCME}).
\end{remark}

\section{Numerical results}\label{sec:7}
In this section, we will demonstrate the robustness of the proposed monolithic mixed-dimensional multigrid method through different numerical experiments. In addition to seeing that the behavior of multigrid is independent of the spatial discretization parameter, we will also analyze the robustness of the algorithm with respect to fracture properties, as the permeability. Further, we will study how the multigrid performance is influenced by considering several fractures, and also illustrate the good behavior of the method on a benchmark problem from the literature. Throughout the whole section, we will consider $\xi=1$ in the model and we will use a $W$-cycle with two pre- and two post-smoothing steps, since this choice has been shown to provide very good results when monolithic multigrid solvers are considered for coupled problems \cite{PeiyaoSISC,PeiyaoJCP}. In our case, we will see that it gives \textit{multigrid textbook efficiency} \cite{TOS01}.

\subsection{One fracture test} 
\begin{figure}[t]
	\begin{center}
		\begin{tabular}{cc}
			\hspace{-0.3cm}\includegraphics[scale = 0.415]{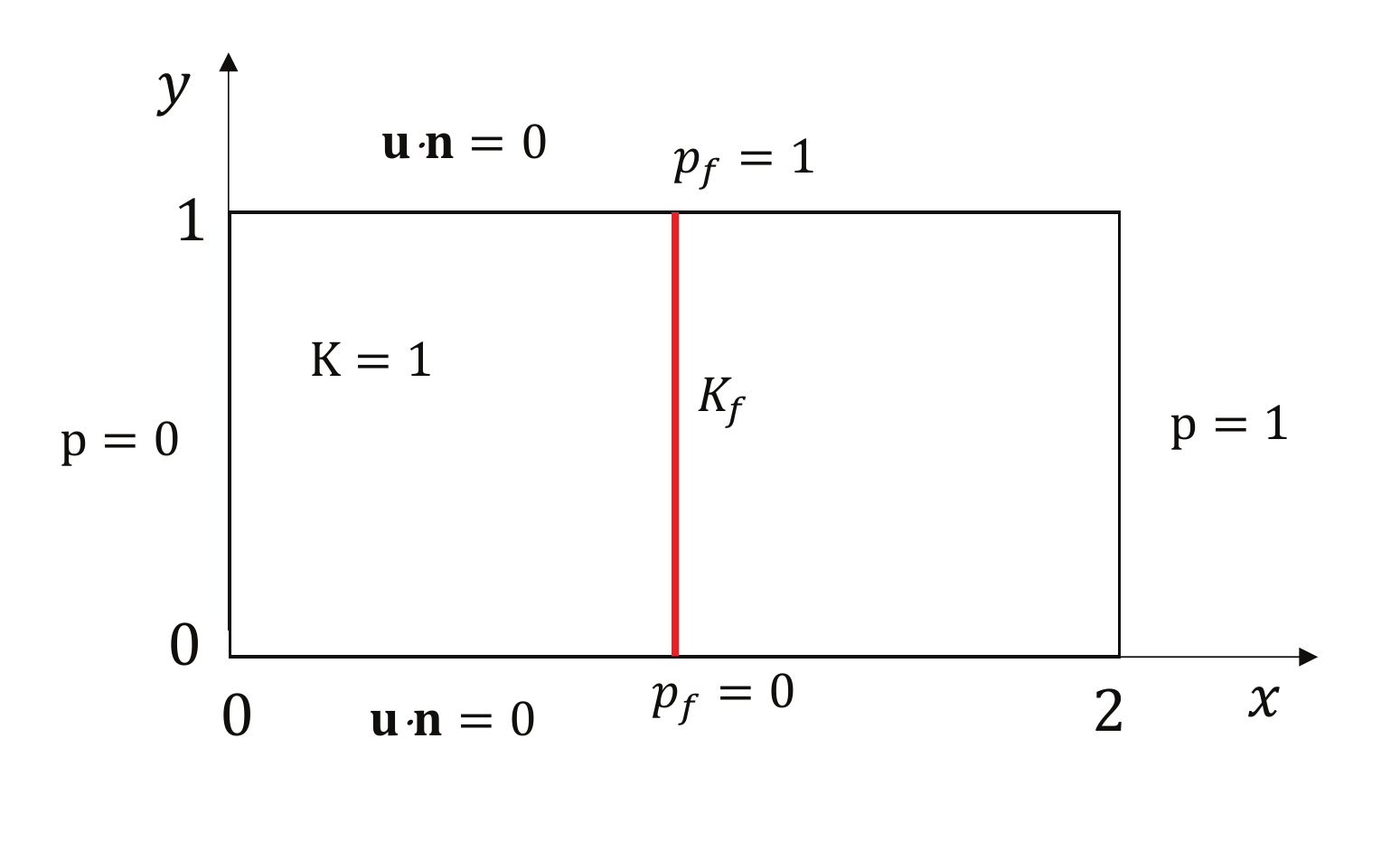}
			&
			\hspace{-0.6cm}\includegraphics[scale = 0.415]{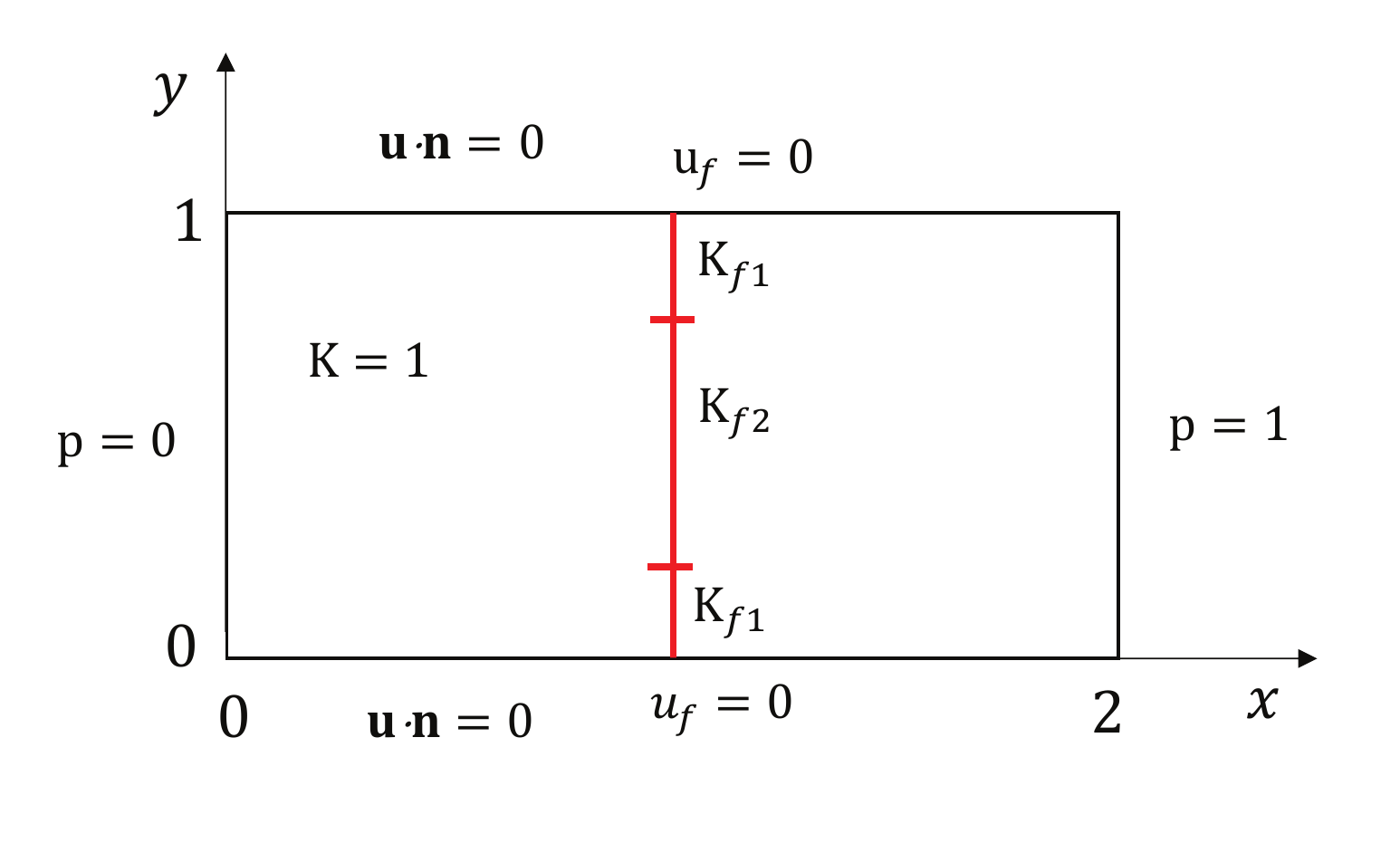}
			\\
			(a) Case 1 & (b) Case 2
		\end{tabular}
	\end{center}
	\label{first_experiment_domain_bc}
	\caption{Domain and boundary conditions for the first numerical experiment.}
\end{figure}

We first consider a test problem presented in \cite{Martin-Jaffre-Roberts} in which the domain consists of an horizontal rectangular slice of porous medium $\Omega = (0,2)\times (0,1)$ with unit permeability (${\mathbf K} = K\mathbf{I}$, with $K=1$ and $\mathbf{I}$ the identity tensor), impermeable bottom and top boundaries and a prescribed pressure of zero and one in the left and right sides, respectively. 
Such domain is divided into two equally sized subdomains by a vertical fracture of width $d=10^{-2}$ for which we consider two different cases: constant permeability in the whole fracture and different values of the permeability within the fracture. Also two different types of boundary conditions are considered in the extremities of the fracture. All these settings are displayed in Figure \ref{first_experiment_domain_bc} for both Case 1 and Case 2.

\subsubsection{Case 1: constant permeability in the fracture}
\begin{figure}[t]
	\begin{center}
		\begin{tabular}{cc}
			\includegraphics[scale = 0.3]{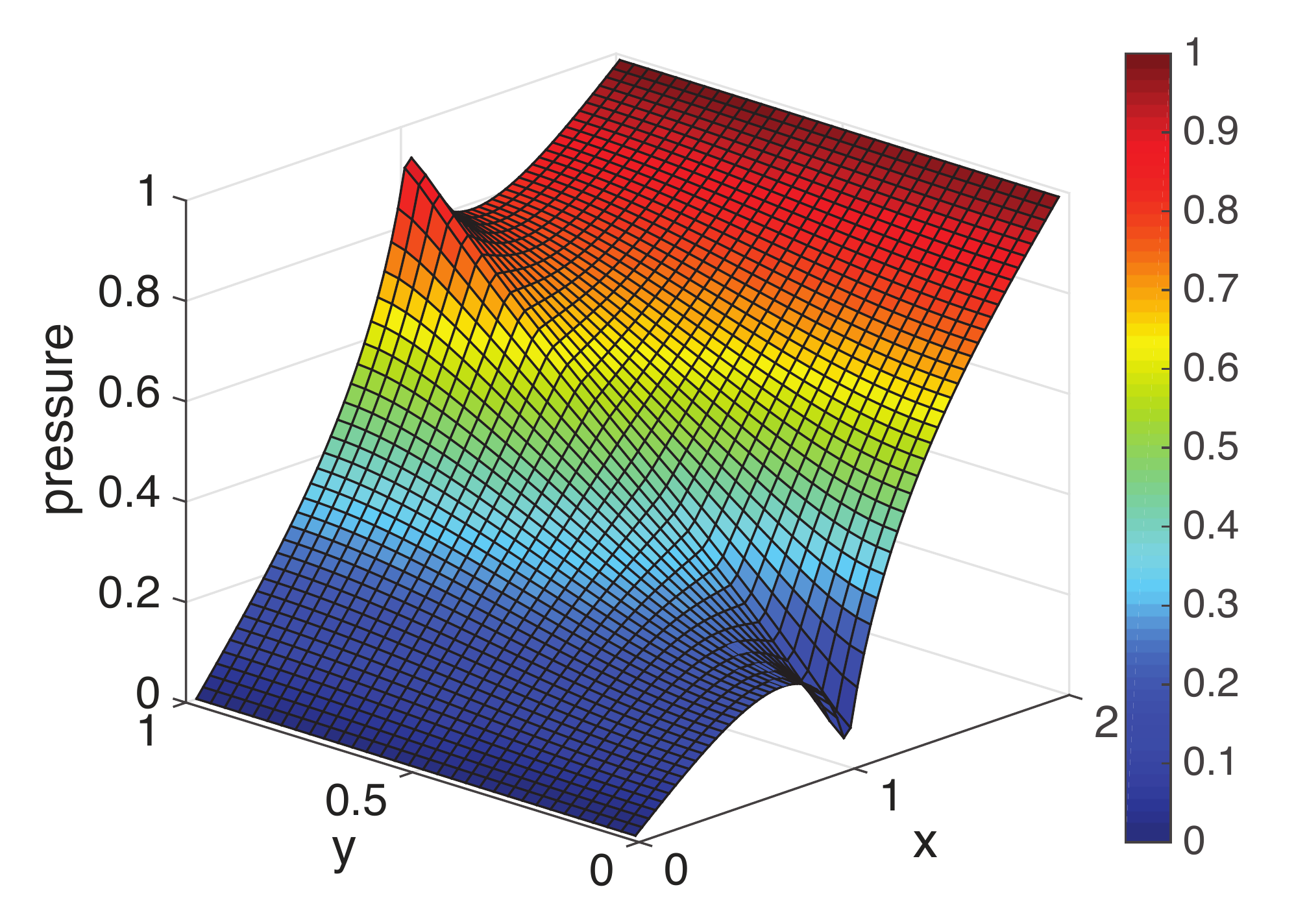}
			&
			\includegraphics[scale = 0.3]{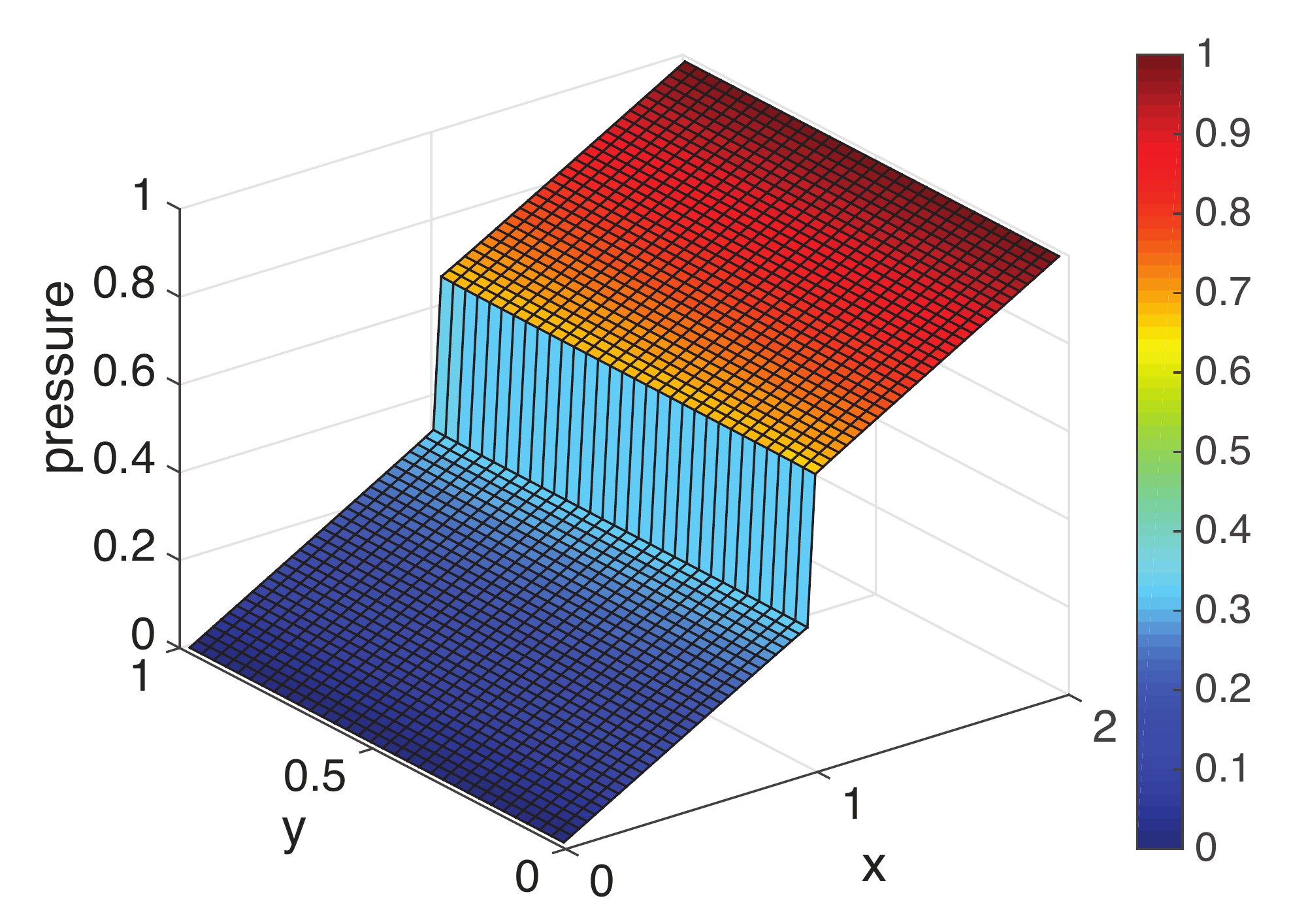}\\
			(a) & (b)
		\end{tabular}
	\end{center}
	\caption{Pressure solution for the fracture permeabilities (a) $K_f = 100$ and (b) $K_f = 0.01$ for the first numerical experiment (Case 1).}
	\label{first_test_pressure}
\end{figure}

In this first case we consider the setting displayed in Figure \ref{first_experiment_domain_bc} (a). The boundary conditions for the fracture are Dirichlet in this case. More concretely, $p_f = 1$ on the top extremity of the fracture, and $p_f = 0$ on the bottom. The permeability tensor in the fracture is given by ${\mathbf K}_f = K_f \mathbf{I}$, and we want to study the influence of different values of $K_f$ on the performance of the multigrid solver. We consider both conductive fractures and blocking fractures, characterized by high or low permeabilities, respectively. As an example, in Figure \ref{first_test_pressure} we show the pressure solution obtained for two different values of $K_f$, one representative of a high permeability (left side) and the other one characteristic of a low permeability (right side) in the fracture. 
Now, to study the robustness of the mixed-dimensional multigrid method with respect to different values of the permeability, in Table \ref{test1_different_permeabilities} we display the number of iterations needed to reduce the initial residual in a factor of $10^{-10}$ for different grid sizes and for low and high permeabilities. 
\begin{table}[t]
	\begin{center}
		\caption{Number of iterations of the mixed-dimensional multigrid method necessary to solve the first numerical experiment (Case 1) with different values of a constant permeability tensor in the fracture.}
		\begin{tabular}{ccccccc}
			\cline{2-7}
			& $K_f$ & $32\times 16$ & $64\times 32$ & $128\times 64$ & $256\times 128$ & $512\times 256$ \\
			\hline
			\multirow{3}{*}{\begin{minipage}{2.3cm}\begin{center} low \\ permeability \end{center} \end{minipage}} & $10^{-6}$ & 8 & 8 & 9 & 9 & 9 \\
			& $10^{-4}$ & 8 & 8 & 9 & 9 & 9 \\ 
			& $10^{-2}$ & 8 & 8 & 9 & 9 & 9 \\
			\hline
			\hline
			\multirow{3}{*}{\begin{minipage}{2.3cm}\begin{center} high \\ permeability \end{center} \end{minipage}} & $10^{2}$ & 10 & 9 & 9 & 10 & 10 \\
			& $10^{4}$ & 8 & 9 & 9 & 9 & 10 \\ 
			& $10^{6}$ & 8 & 9 & 9 & 9 & 10 \\
			\hline
		\end{tabular}
	\end{center}
	\label{test1_different_permeabilities}
\end{table}
We can observe that, for all the values of $K_f$, the performance of the multigrid method is independent of the spatial discretization parameter. Moreover, only a few iterations are required to satisfy the stopping criterion.

\subsubsection{Case 2: variable permeability in the fracture}
\begin{figure}[t]
	\begin{center}
		\begin{tabular}{cc}
			\includegraphics[scale = 0.3]{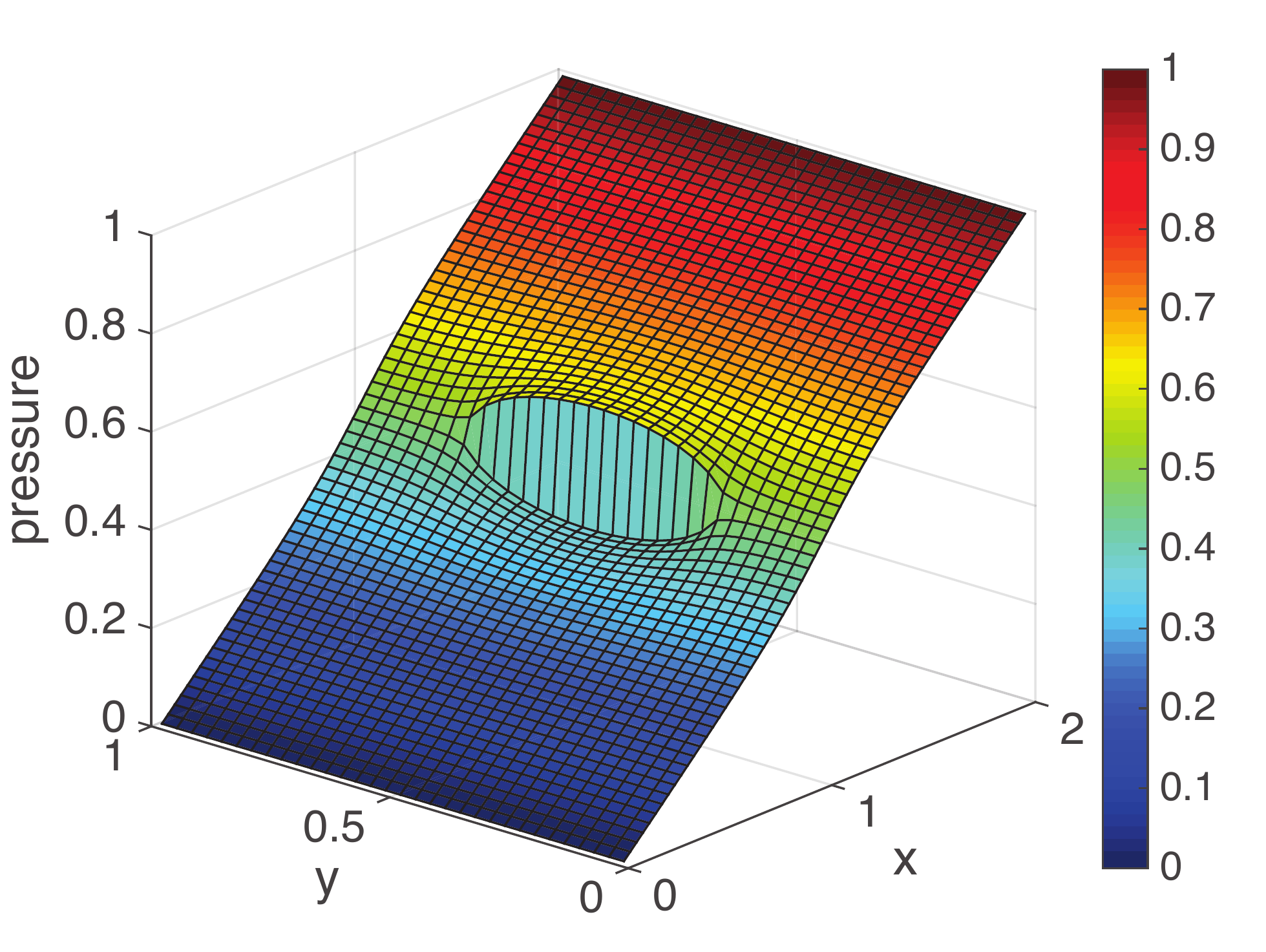}
			&
			\includegraphics[scale = 0.31]{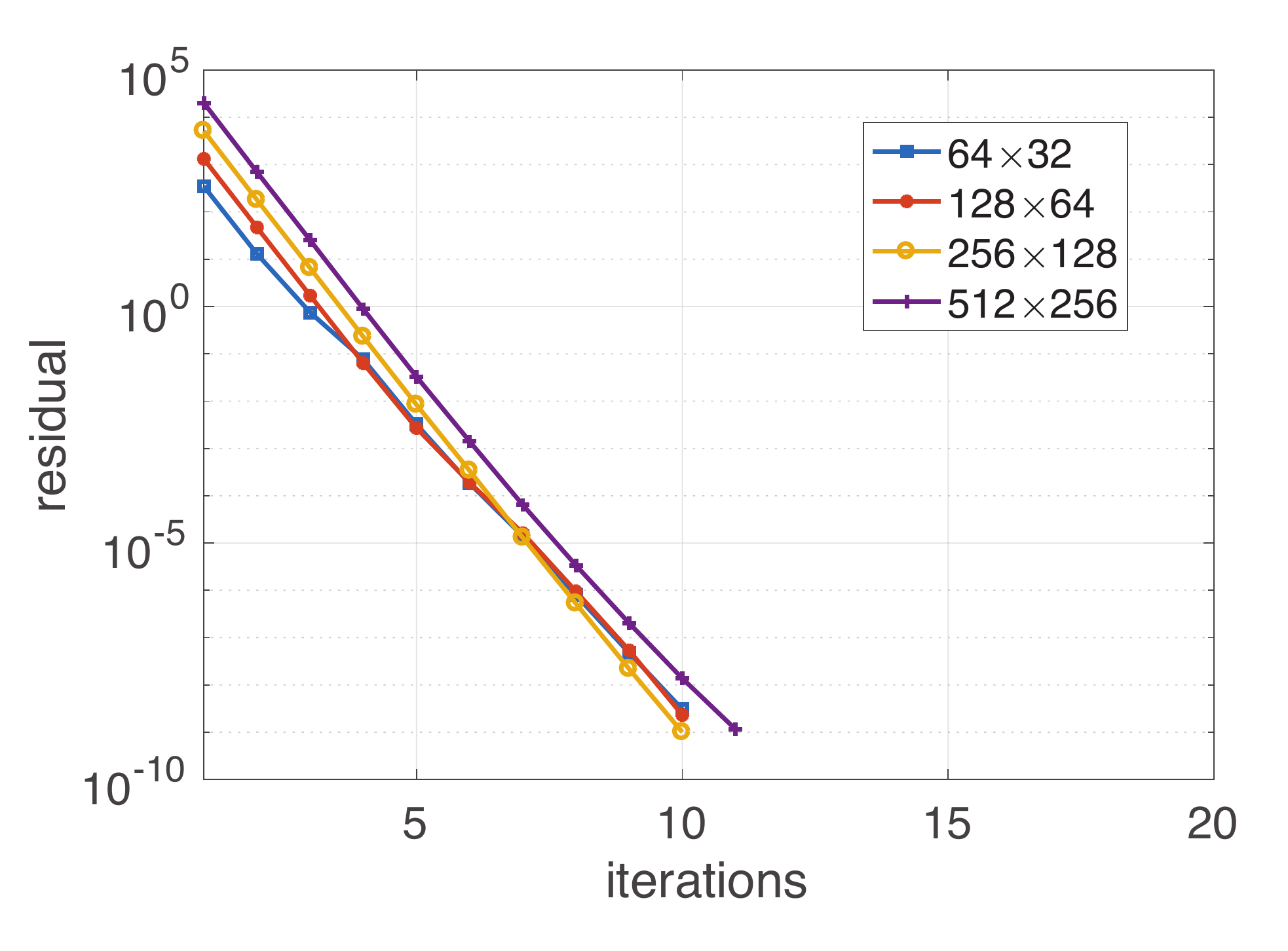}\\
			(a) & (b)
		\end{tabular}
	\end{center}
	\caption{(a) Pressure solution and (b) history of the convergence of the proposed multigrid method for the first numerical experiment (Case 2).}
	\label{first_test_case2}
\end{figure}

We consider now the setting displayed in Figure \ref{first_experiment_domain_bc} (b). In this case, the boundary conditions for the fracture are homogeneous Neumann conditions on both extremities of the fracture. The permeability tensor in the fracture is now given by 
$$
{\mathbf K}_f =\begin{cases}
K_{f1} \mathbf{I},& \hbox{if } \; 0<y<\frac{1}{4} \; \hbox{ or } \; \frac{3}{4}<y<1,\\[1ex]
K_{f2} \mathbf{I},& \hbox{if } \; \frac{1}{4}<y<\frac{3}{4}.
\end{cases}
$$
In particular, we consider $K_{f1} = 10^2$ and $K_{f2} = 2\times10^{-3}$. Due to the low value of $K_{f2}$ the fluid tends to avoid the middle part of the fracture, representing a geological barrier. This behavior can be clearly seen in the pressure distribution depicted in Figure \ref{first_test_case2} (a). 
Finally, we want to study if this changes of permeability within the fracture have some effect on the multigrid performance. For this purpose, in Figure \ref{first_test_case2} (b) we display the history of the convergence of the multigrid solver for different mesh sizes. More concretely, the reduction of the residual is depicted against the number of iterations, and the stopping criterion is to reduce the initial residual until $10^{-8}$. It is clearly seen that the convergence of the monolithic mixed-dimensional multigrid method is independent of the spatial discretization parameter, and the number of iterations are very similar to those obtained in the previous case with a constant permeability in the fracture.

\subsection{Four fracture network}
In this numerical experiment we want to see the performance of the proposed multigrid method when several fractures are considered. For this purpose, we consider a network of four fractures whose width is $d=10^{-2}$. We perform two different tests corresponding to non-connected and connected fractures. In both cases, we consider a unit square porous medium domain with impermeable lateral walls and a given pressure on top ($p=1$) and bottom ($p=0$). The permeability of the porous matrix is given by the identity tensor and the permeability in the fractures is given by $\mathbf{K}_{fi}=K_{fi}\mathbf{I}$, with $K_{fi}$ certain constants to be determined below, for $i=1,2,3,4$. Following \cite{ang:boy:hub:09,sco:for:sot:17}, immersed fracture tips are modeled by homogeneous flux conditions.

\subsubsection{Case 1: four non-connected fractures}
\begin{figure}[t]
	\begin{center}
		\begin{tabular}{cc}
			\begin{minipage}{0.5\textwidth}
				\begin{center}
					\includegraphics[scale = 0.43]{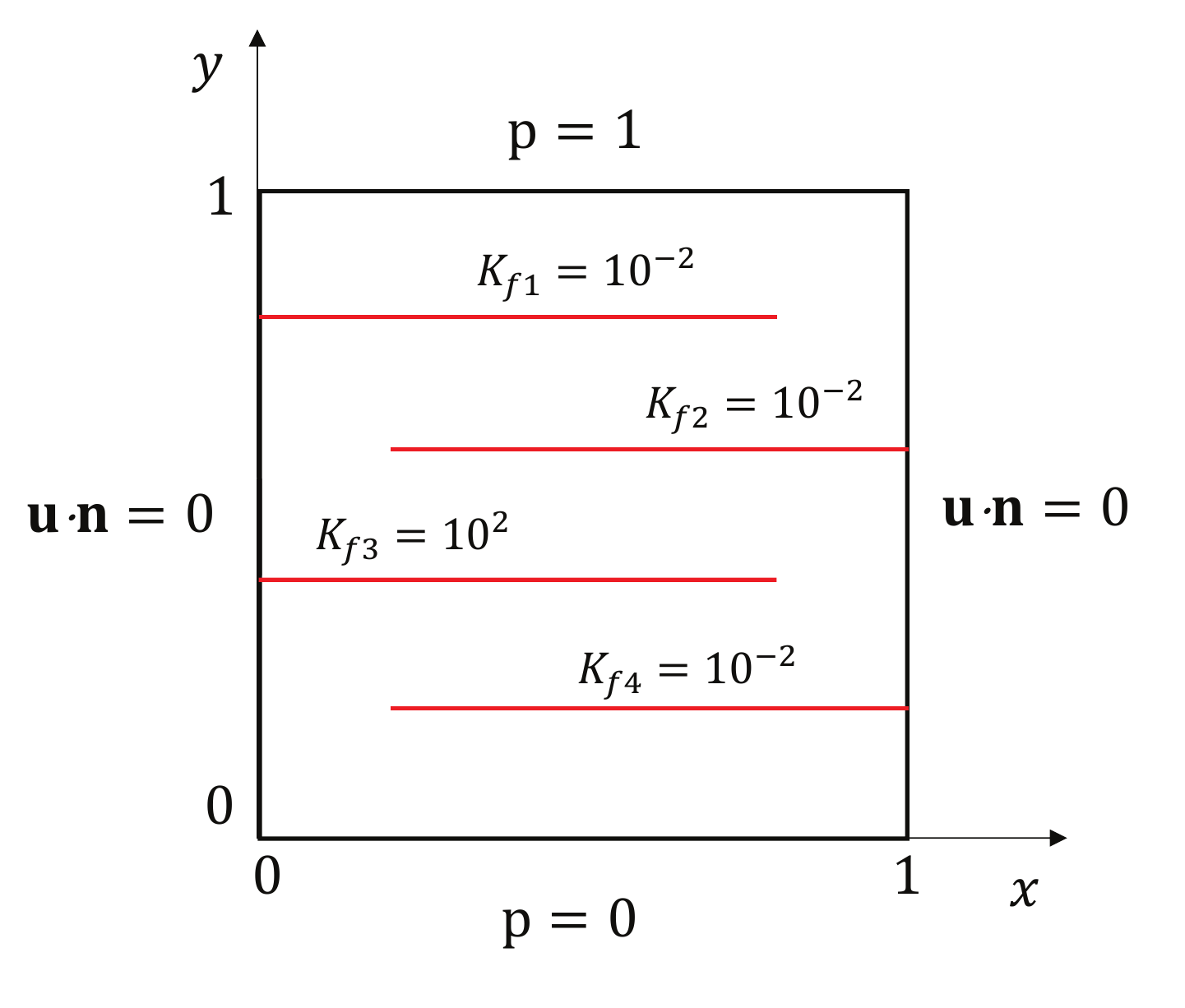}
				\end{center}
			\end{minipage}
			&
			\begin{minipage}{0.43\textwidth}
				\begin{center}
					\includegraphics[scale = 0.28]{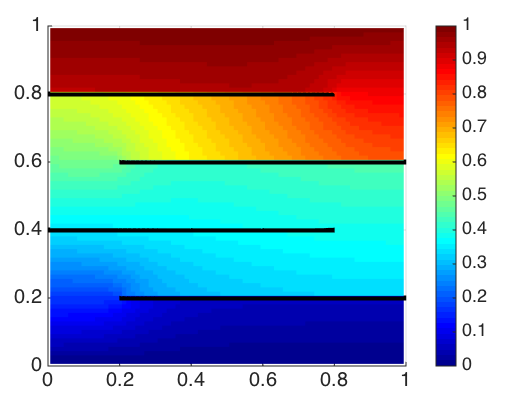}
				\end{center}
			\end{minipage}
			\\
			(a) & (b)
		\end{tabular}
	\end{center}
	\caption{(a) Fracture network and settings, and (b) corresponding pressure distribution for the four fracture experiment (Case 1).}
	\label{four_fractures}
\end{figure}

The settings considered in this first case are based on \cite{ang:boy:hub:09} and schematized in Figure \ref{four_fractures} (a). A set of four horizontal fractures $\left\{{\gamma}_i\right\}_{i=1,\ldots,4}$ is considered, where:
\begin{eqnarray*}
	{\gamma}_1 &=& \left\{ (x,y) \; | \; y = 0.8,\; 0\leq x\leq 0.8 \right\}, \\
	{\gamma}_2 &=& \left\{ (x,y) \; | \; y = 0.6,\; 0.2\leq x\leq 1 \right\}, \\
	{\gamma}_3 &=& \left\{ (x,y) \; | \; y = 0.4,\; 0\leq x\leq 0.8 \right\}, \\
	{\gamma}_4 &=& \left\{ (x,y) \; | \; y = 0.2,\; 0.2\leq x\leq 1 \right\}. 
\end{eqnarray*} 
Fractures $\gamma_1$, $\gamma_2$ and $\gamma_4$ are assumed to be barriers with $K_{f1} = K_{f2}=K_{f4}=10^{-2}$, whereas fracture $\gamma_3$ is highly conductive with $K_{f3}=10^2$. 
The effect of these fractures on the pressure distribution within the porous medium domain can be seen in Figure \ref{four_fractures} (b).

\subsubsection{Case 2: four connected fractures}
\begin{figure}[t]
	\begin{center}
		\begin{tabular}{cc}
			\begin{minipage}{0.5\textwidth}
				\begin{center}
					\includegraphics[scale = 0.43]{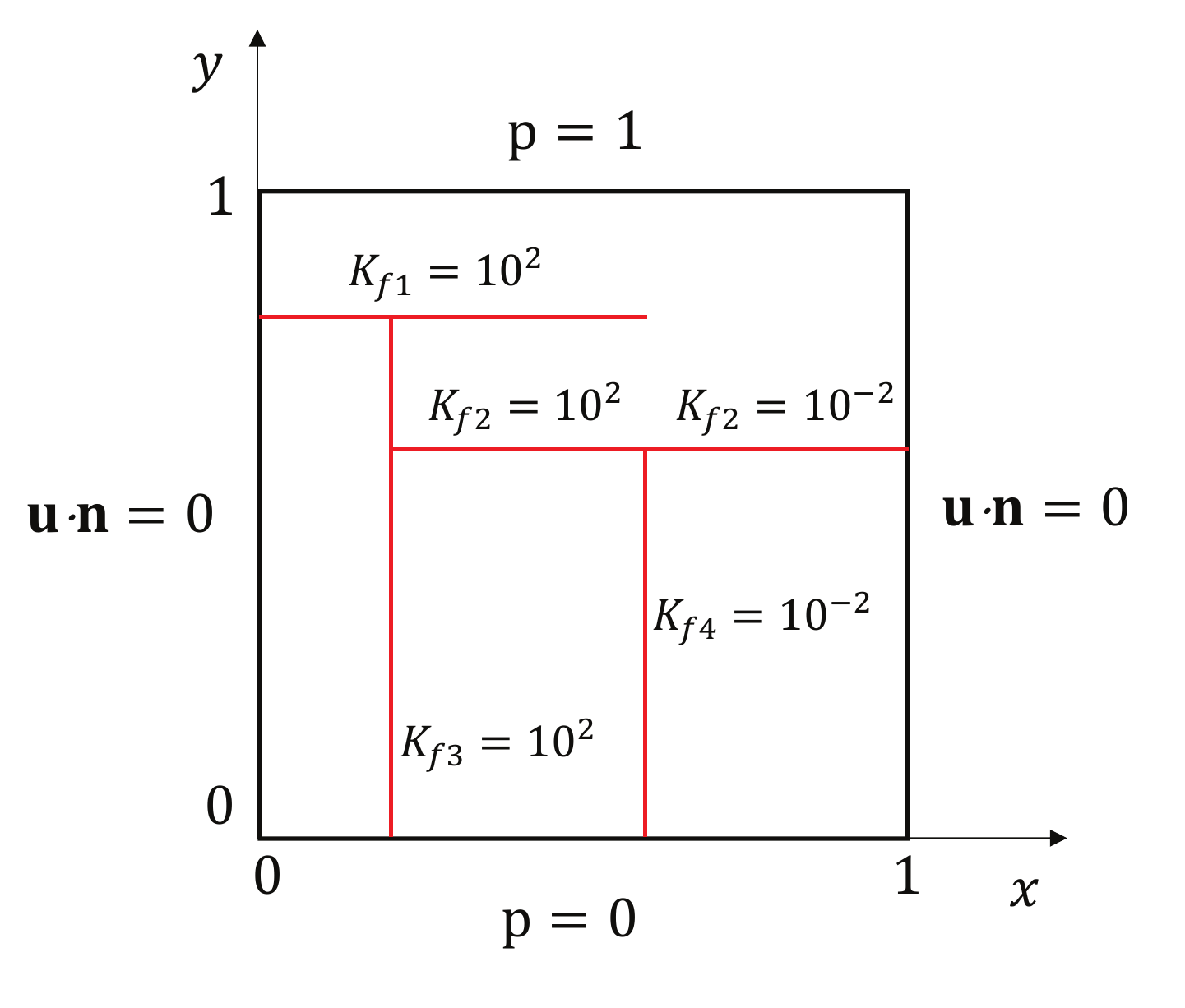}
				\end{center}
			\end{minipage}
			&
			\begin{minipage}{0.43\textwidth}
				\begin{center}
					\includegraphics[scale = 0.28]{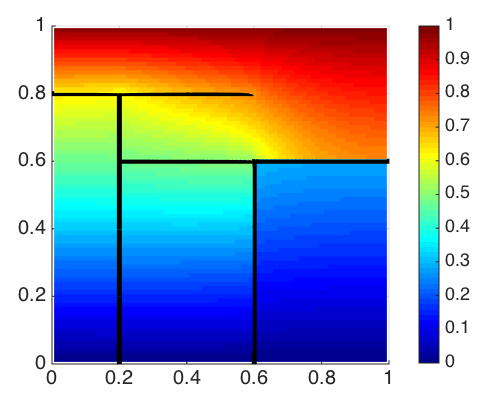}
				\end{center}
			\end{minipage}
			\\
			(a) & (b)
		\end{tabular}
	\end{center}
	\caption{(a) Fracture network and settings, and (b) corresponding pressure distribution for the four fracture experiment (Case 2).}
	\label{four_fractures_cross}
\end{figure}

In the second case, we assume a fracture network with four fractures which are connected. A schematic picture of the network together with the properties of the fractures can be seen in Figure \ref{four_fractures_cross} (a). 
In particular, we consider two vertical fractures 
\begin{eqnarray*}
	\gamma_3 &=& \left\{ (x,y) \; | \; x = 0.2, \; 0\leq y \leq 0.8\right\},\\
	\gamma_4 &=& \left\{ (x,y) \; | \; x = 0.6, \; 0\leq y \leq 0.6\right\},
\end{eqnarray*}
with constant permeabilities given by $K_{f3}=10^2$ and $K_{f4} = 10^{-2}$, respectively. Moreover, we consider two horizontal fractures 
\begin{eqnarray*}
	\gamma_1 &=& \left\{ (x,y) \; | \; y =0.8, \; 0\leq x \leq 0.6\right\},\\
	\gamma_2 &=& \left\{ (x,y) \; | \; y = 0.6, \; 0.2\leq x \leq 1\right\}.
\end{eqnarray*}
The first one has a constant permeability given by $K_{f1} = 10^2$ whereas the second one has a variable permeability given by:
$$
K_{f2} =
\begin{cases}
10^2,& \hbox{if } \; 0.2\leq x \leq 0.6,\\
10^{-2},& \hbox{if } \; 0.6<x\leq 1.
\end{cases}
$$
The pressure solution is shown in Figure \ref{four_fractures_cross} (b), where we can observe the effect of the different permeabilities of the fractures on the pressure at the porous matrix. 

We perform the proposed monolithic multigrid method for solving both test cases. Table \ref{test_four_fractures} shows the number of multigrid iterations needed to reduce the initial residual in a factor of $10^{-10}$  for different mesh sizes. We can see a very robust behavior of the multigrid algorithm for both cases since few iterations are enough to satisfy the stopping criterion. 
\begin{table}[t]
	\begin{center}\caption{Number of iterations of the mixed-dimensional multigrid method necessary to solve the four fracture experiments for different grid sizes.}
		\begin{tabular}{ccccccc}
			\cline{2-7}
			& $40\times 40$ & $80\times 80$ & $160\times 160$ & $320\times 320$ & $640\times 640$ & $1280\times 1280$  \\
			\hline
			Case 1 & 9 & 9 & 10 & 10 & 11 & 11 \\
			\hline
			\hline
			Case 2 & 11 & 11 & 11 & 12 & 13 & 13 \\
			\hline
		\end{tabular}
	\end{center}
	\label{test_four_fractures}
\end{table}
From this experiment, we observe that the monolithic mixed-dimensional multigrid method is also robust when several fractures (connected and/or non-connected) are considered. 

\subsection{Benchmark problem}\label{subsec:benchmark}
\begin{figure}[t]
	\begin{center}
		\centering\includegraphics[scale = 0.5]{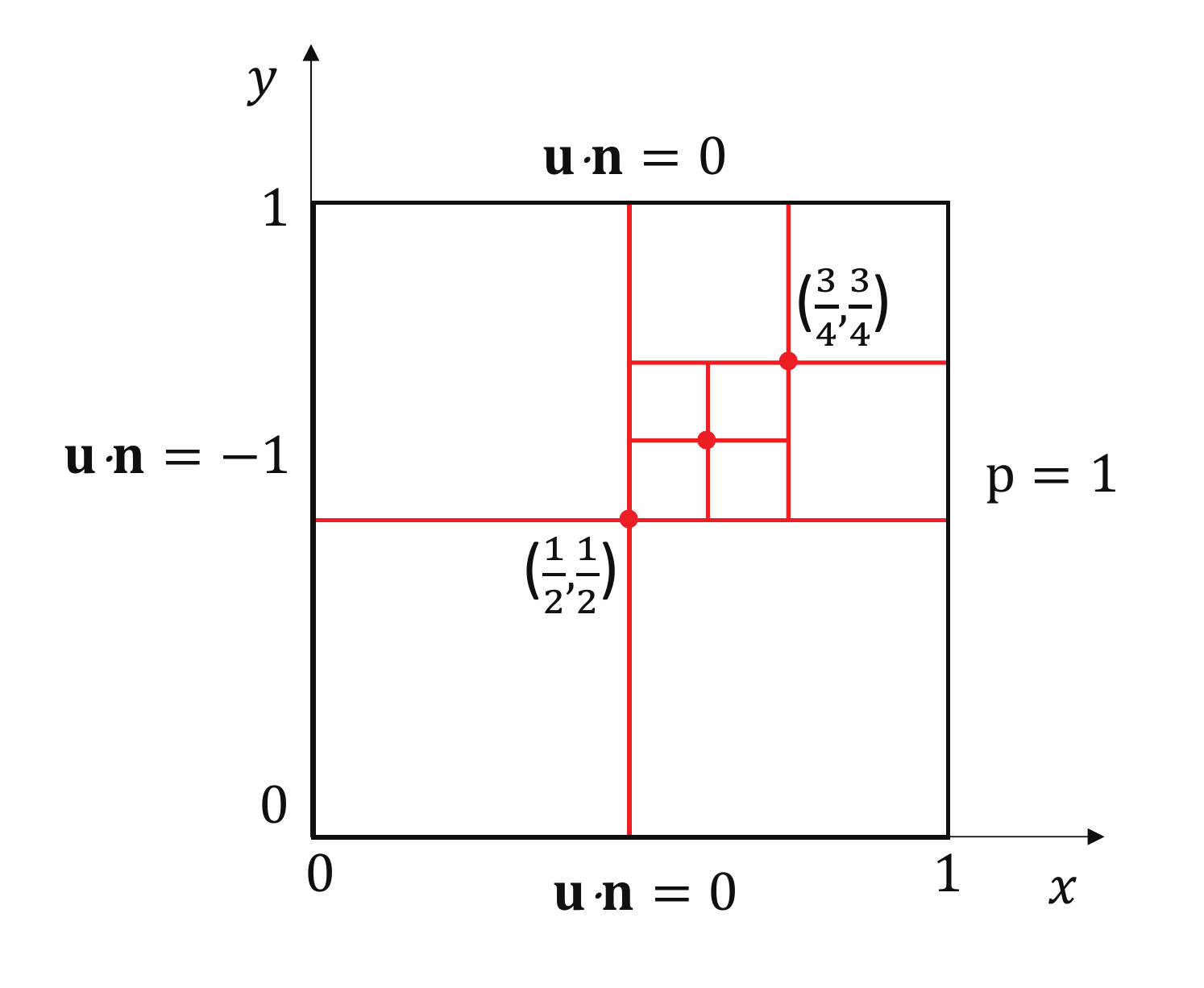}
	\end{center}
	\caption{Configuration and boundary conditions for the network of fractures in the benchmark problem.}
	\label{Figure_benchmark}
\end{figure}

The last numerical experiment considered in this work is a benchmark problem for single-phase flow in fractured porous media stated in \cite{FLEMISCH2018239}. This test is based on a problem proposed in 
\cite{Geiger} with different boundary conditions and material properties. The fracture network embedded in the unit square domain is shown in Figure \ref{Figure_benchmark}. 
We have the Dirichlet condition $p = 1$ on the right side. Homogeneous Neumann boundary conditions on the top and bottom are prescribed, whereas on the left side, we have a constant inflow flow, ${\mathbf u} \cdot {\mathbf n} = -1$ (see Figure \ref{Figure_benchmark}). 
The permeability matrix is fixed to ${\bf K} = {\bf I}$, and all the fractures have a constant width $d = 10^{-4}$.
As in \cite{FLEMISCH2018239}, we consider two cases for the permeability tensor in the fracture ${\mathbf K}_f = K_f \mathbf{I}$: a case 
where the fluid tends to flow rapidly along the fracture, i.e., $K_{f}=10^4$; and a case where the fluid tends to avoid the fracture, i.e., $K_{f}=10^{-4}$. Both cases have been solved by the proposed mixed-dimensional multigrid method, giving rise to a very efficient solver independent of the fracture network. 

\begin{figure}[t]
	\begin{center}
		\begin{tabular}{cc}
			\hspace{-0.4cm}\includegraphics[scale = 0.17]{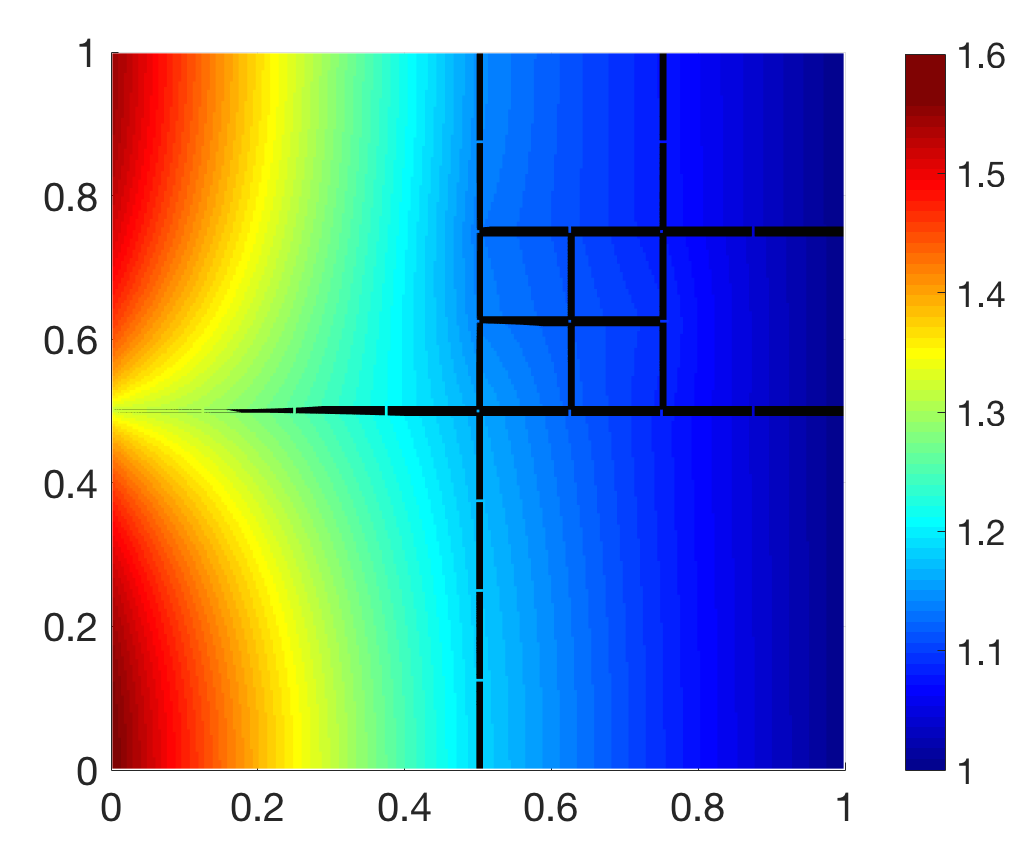}
			&
			\includegraphics[scale = 0.346]{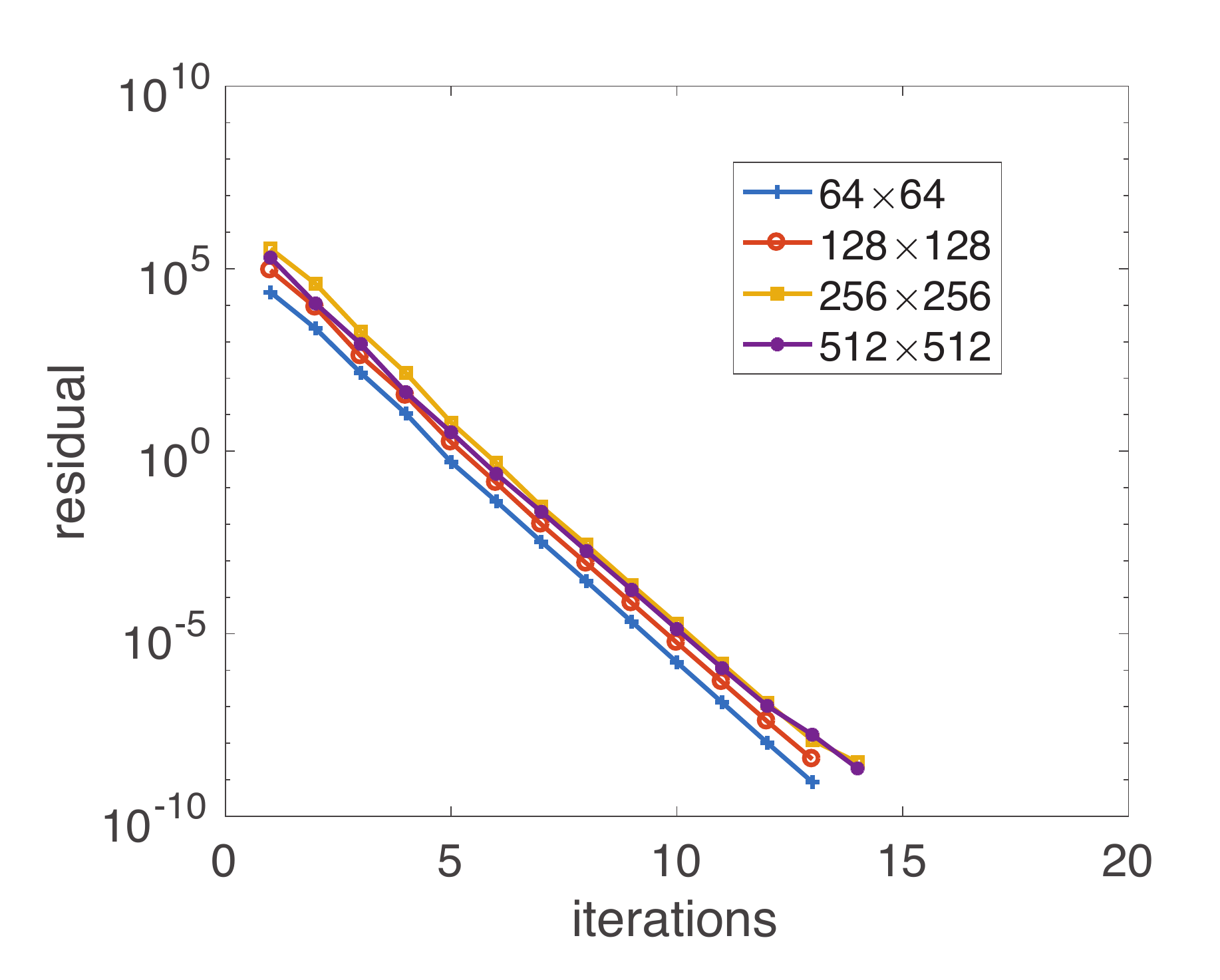}\\
			(a) & (b)
		\end{tabular}
	\end{center}
	\caption{(a) Pressure solution and (b) history of the convergence of the proposed multigrid method for the benchmark problem with highly conductive fractures.}
	\label{benchmark_case1}
\end{figure}

For the first case in which we deal with a highly conductive network, the obtained pressure distribution is displayed in Figure \ref{benchmark_case1} (a). It is clear its correspondence with the reference solution of the benchmark problem published in \cite{FLEMISCH2018239}. 
The history of the convergence of the monolithic multigrid method is depicted in Figure \ref{benchmark_case1} (b), where the residual reduction is shown for different grid sizes. The stopping criterion is to reduce the initial residual until $10^{-8}$. We can observe that, as expected, the performance of the multigrid method is independent of the spatial discretization parameter. Moreover, it results in a very efficient solver since only around $10$ iterations are enough to solve the whole fracture network. 

\begin{figure}
	\begin{center}
		\begin{tabular}{cc}
			\hspace{-0.4cm}\includegraphics[scale = 0.17]{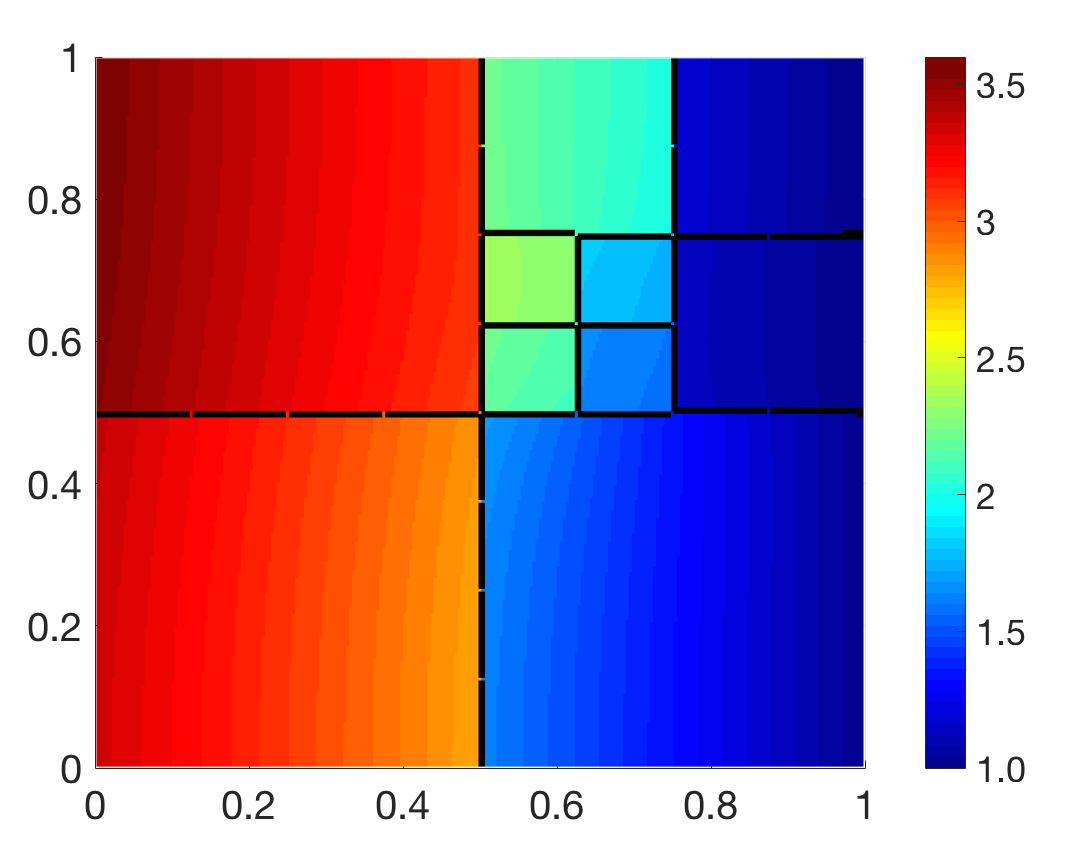}
			&
			\hspace{-0.3cm}\includegraphics[scale = 0.34]{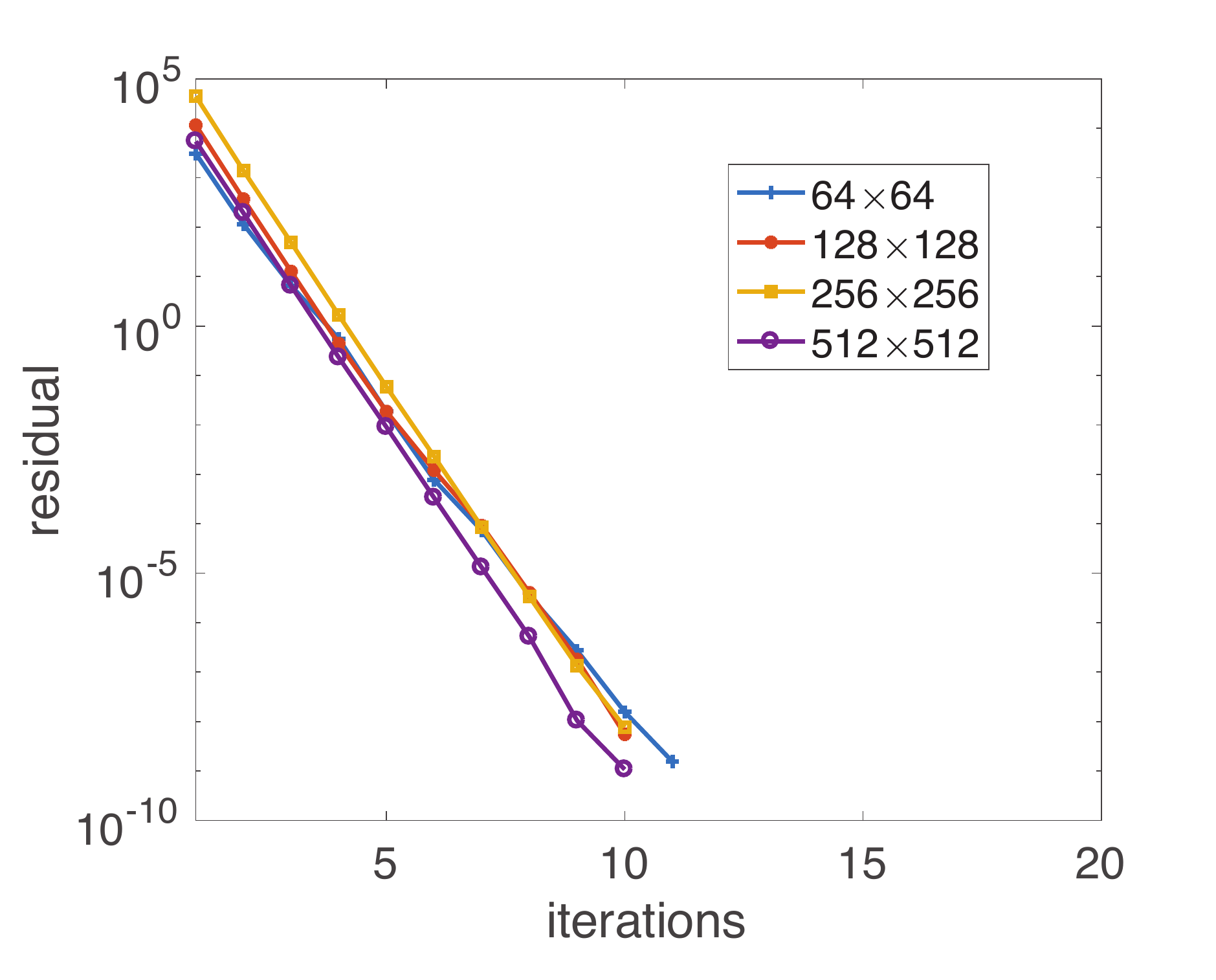}\\
			(a) & (b)
		\end{tabular}
	\end{center}
	\caption{(a) Pressure solution and (b) history of the convergence of the proposed multigrid method for the benchmark problem with blocking fractures.}
	\label{benchmark_case2}
\end{figure}

Similar results are obtained for the second case, corresponding to a blocking fracture network. Again, the pressure distribution matches perfectly with the reference solution in \cite{FLEMISCH2018239}, as shown in Figure \ref{benchmark_case2} (a). We can observe in the picture the pressure discontinuities reminiscent of the low permeability in the fractures. 
In Figure \ref{benchmark_case2} (b), we display the history of the convergence of the proposed mixed-dimensional multigrid method for this second case of the benchmark problem. 
The convergence results for the blocking fracture network are similar to those presented for the conducting fractures. For both cases, the monolithic mixed-dimensional multigrid method shows an excellent performance.

\begin{remark}
	Local Fourier analysis \cite{Bra77, Bra94,Wie01} is the main quantitative analysis to predict the convergence rates of multigrid algorithms. 
	LFA assumes that all operators involved in the multigrid procedure are local, have constant coefficients and are defined on an infinite grid neglecting the effect of boundary conditions.
	It seems not possible, or at least it is not clear, how to carry out a local Fourier analysis for the mixed-dimensional multigrid method proposed here. Nevertheless, we have performed LFA to predict the convergence of the multigrid based on the Vanka smoother for the Darcy problem considered in the porous matrix. The analysis of Vanka-type smoothers by LFA requires a special strategy and here we have followed the ideas presented in \cite{LFA_overlapping}. 
	From this analysis, we observe that the convergence rates obtained for the whole problem including fractures are very close to those provided in the case of simple Darcy flow. More concretely, by using four smoothing steps, a two-grid convergence factor of $0.04$ is obtained for this latter case, whereas the asymptotic rates observed in the numerical experiments carried out in this work vary from $0.04$ in the first experiment to $0.085$ in the benchmark problem which is the worst case. This means that the implementation and the treatment in the fractures is done in an optimal way, since the convergence of the whole fracture network problem is very similar to that for the Darcy flow problem.  
\end{remark} 

\section{Conclusions}\label{sec:8}
We have proposed a monolithic mixed-dimensional multigrid method for solving single-phase flow problems in porous media with intersecting fractures. The algorithm is based on combining two-dimensional smoother and inter-grid transfer operators in the porous matrix with their one-dimensional counterparts within the fracture network. This exotic union gives rise to a very efficient solver for this type of problems. The robustness of the proposed method with respect to different parameters of the fractures, as the permeability, as well as to the number of fractures and the size of the grid has been shown through different numerical experiments, including a benchmark problem from the literature.

\bibliographystyle{plain}
\bibliography{references}

\begin{thebibliography}{10}

\bibitem{ahm:edw:lam:hui:pal:15a}
R.~Ahmed, M.~G. Edwards, S.~Lamine, B.~A.~H. Huisman, and M.~Pal.
\newblock Control-volume distributed multi-point flux approximation coupled
  with a lower-dimensional fracture model.
\newblock {\em J. Comput. Phys.}, 284:462--489, 2015.

\bibitem{alb:jaf:rob:ser:01}
Clarisse Alboin, J\'er\^ome Jaffr\'e, Jean~E. Roberts, and Christophe Serres.
\newblock Modeling fractures as interfaces for flow and transport in porous
  media.
\newblock In {\em Fluid Flow and Transport in Porous Media: Mathematical and
  Numerical Treatment}, volume 295 of {\em Contemp. Math.}, pages 13--24. Amer.
  Math. Soc., Providence, RI, 2002.

\bibitem{ami:ker:mar:rob:2006}
Laila Amir, Michel Kern, Jean~E Roberts, and Vincent Martin.
\newblock D{\'e}composition de domaine pour un milieu poreux fractur{\'e}: un
  mod{\`e}le en 3{D} avec fractures qui s'intersectent.
\newblock {\em ARIMA}, 5:11--25, 2006.

\bibitem{ang:boy:hub:09}
Philippe Angot, Franck Boyer, and Florence Hubert.
\newblock Asymptotic and numerical modelling of flows in fractured porous
  media.
\newblock {\em M2AN Math. Model. Numer. Anal.}, 43(2):239--275, 2009.

\bibitem{ant:fac:rus:ver:16b}
Paola~F. Antonietti, C.~Facciol\`a, A.~Russo, and Marco Verani.
\newblock Discontinuous {G}alerkin approximation of flows in fractured porous
  media on polytopic grids.
\newblock Technical Report MOX-55/2016, Dipartimento di Matematica, Politecnico
  di Milano, 2016.

\bibitem{ant:for:sco:ver:ver:16}
Paola~F. Antonietti, Luca Formaggia, Anna Scotti, Marco Verani, and Nicola
  Verzott.
\newblock Mimetic finite difference approximation of flows in fractured porous
  media.
\newblock {\em ESAIM Math. Model. Numer. Anal.}, 50(3):809--832, 2016.

\bibitem{arb:dou:90}
Todd Arbogast and Jim Douglas, Jr.
\newblock Dual-porosity models for flow in naturally fractured reservoirs.
\newblock In J.~Cushman, editor, {\em Dynamics of Fluids in Hierarchical Porous
  Media}, pages 177--221. Academic Press, London, 1990.

\bibitem{arb:dou:hor:90}
Todd Arbogast, Jim Douglas, Jr., and Ulrich Hornung.
\newblock Derivation of the double porosity model of single phase flow via
  homogenization theory.
\newblock {\em SIAM J. Math. Anal.}, 21(4):823--836, 1990.

\bibitem{arb:whe:yot:97}
Todd Arbogast, Mary~F Wheeler, and Ivan Yotov.
\newblock Mixed finite elements for elliptic problems with tensor coefficients
  as cell-centered finite differences.
\newblock {\em SIAM J. Numer. Anal.}, 34(2):828--852, 1997.

\bibitem{ben:ber:pie:sci:14}
Mat\'ias~Fernando Benedetto, Stefano Berrone, Sandra Pieraccini, and Stefano
  Scial\`o.
\newblock The virtual element method for discrete fracture network simulations.
\newblock {\em Comput. Methods Appl. Mech. Engrg.}, 280:135--156, 2014.

\bibitem{boo:nor:vat:17}
Wietse~M. Boon, Jan~M. Nordbotten, and Jon~E. Vatne.
\newblock Mixed-dimensional elliptic partial differential equations.
\newblock {\em arXiv:1710.00556 [math.AP]}, 2017.

\bibitem{boo:nor:yot:18}
Wietse~M. Boon, Jan~M. Nordbotten, and Ivan Yotov.
\newblock Robust discretization of flow in fractured porous media.
\newblock {\em SIAM J. Numer. Anal.}, 56:2203--2233, 2018.

\bibitem{Bra77}
Achi Brandt.
\newblock Multi-level adaptive solutions to boundary-value problems.
\newblock {\em Math. Comp.}, 31(138):333--390, 1977.

\bibitem{Bra94}
Achi Brandt.
\newblock Rigorous quantitative analysis of multigrid. {I}. constant
  coefficients two-level cycle with ${L}_2$-norm.
\newblock {\em SIAM J. Numer. Anal.}, 31(6):1695--1730, 1994.

\bibitem{buk:yot:zun:17}
Martina Buka\v{c}, Ivan Yotov, and Paolo Zunino.
\newblock Dimensional model reduction for flow through fractures in poroelastic
  media.
\newblock {\em ESAIM Math. Model. Numer. Anal.}, 51(4):1429--1471, 2017.

\bibitem{cha:pie:for:18}
Florent Chave, Daniele~A. Di~Pietro, and Luca Formaggia.
\newblock A hybrid high-order method for {D}arcy flows in fractured porous
  media.
\newblock {\em SIAM J. Sci. Comput.}, 40(2):A1063--A1094, 2018.

\bibitem{ang:sco:12}
Carlo D'Angelo and Anna Scotti.
\newblock A mixed finite element method for {D}arcy flow in fractured porous
  media with non-matching grids.
\newblock {\em ESAIM Math. Model. Numer. Anal.}, 46:465--489, 2012.

\bibitem{delpra:fum:sco:17}
Marco Del~Pra, Alessio Fumagalli, and Anna Scotti.
\newblock Well posedness of fully coupled fracture/bulk {D}arcy flow with
  {XFEM}.
\newblock {\em SIAM J. Numer. Anal.}, 55(2):785--811, 2017.

\bibitem{eym:gal:gui:her:mas:14}
R.~Eymard, T.~Gallou\"et, C.~Guichard, R.~Herbin, and R.~Masson.
\newblock T{P} or not {TP}, that is the question.
\newblock {\em Comput. Geosci.}, 18(3-4):285--296, 2014.

\bibitem{FLEMISCH2018239}
Bernd Flemisch, Inga Berre, Wietse Boon, Alessio Fumagalli, Nicolas Schwenck,
  Anna Scotti, Ivar Stefansson, and Alexandru Tatomir.
\newblock Benchmarks for single-phase flow in fractured porous media.
\newblock {\em Adv. Water Resour.}, 111:239--258, 2018.

\bibitem{fle:fum:sco:16}
Bernd Flemisch, Alessio Fumagalli, and Anna Scotti.
\newblock A review of the {XFEM}-based approximation of flow in fractured
  porous media.
\newblock In {\em Advances in Discretization Methods}, volume~12 of {\em SEMA
  SIMAI Springer Ser.}, pages 47--76. Springer, Cham, 2016.

\bibitem{for:fum:sco:ruf:14}
Luca Formaggia, Alessio Fumagalli, Anna Scotti, and Paolo Ruffo.
\newblock A reduced model for {D}arcy's problem in networks of fractures.
\newblock {\em ESAIM Math. Model. Numer. Anal.}, 48:1089--1116, 2014.

\bibitem{sco:for:sot:17}
Luca Formaggia, Anna Scotti, and Federica Sottocasa.
\newblock Analysis of a mimetic finite difference approximation of flows in
  fractured porous media.
\newblock {\em ESAIM Math. Model. Numer. Anal.}, 52:595--630, 2018.

\bibitem{fri:mar:rob:saa:12}
Najla Frih, Vincent Martin, Jean~E. Roberts, and Ali Sa\^{a}da.
\newblock Modeling fractures as interfaces with nonmatching grids.
\newblock {\em Comput. Geosci.}, 16:1043--1060, 2012.

\bibitem{fri:rob:saa:08}
Najla Frih, Jean~E. Roberts, and Ali Saada.
\newblock Modeling fractures as interfaces: a model for {F}orchheimer
  fractures.
\newblock {\em Comput. Geosci.}, 12(1):91--104, 2008.

\bibitem{fum:kei:18}
Alessio Fumagalli and Eirik Keilegavlen.
\newblock Dual virtual element method for discrete fractures networks.
\newblock {\em SIAM J. Sci. Comput.}, 40:B228--B258, 2018.

\bibitem{fum:sco:13}
Alessio Fumagalli and Anna Scotti.
\newblock A numerical method for two-phase flow in fractured porous media with
  non-matching grids.
\newblock {\em Adv. Water Resour.}, 62:454--464, 2013.

\bibitem{Geiger}
Sebastian Geiger, M.~Dentz, and I.~Neuweiler.
\newblock A novel multi-rate dual-porosity model for improved simulation of
  fractured and multi-porosity reservoirs.
\newblock {\em SPE J.}, 18(4):670--684, 8 2013.

\bibitem{gir:whe:gan:mea:15}
V.~Girault, M.~F. Wheeler, B.~Ganis, and M.~E. Mear.
\newblock A lubrication fracture model in a poro-elastic medium.
\newblock {\em Math. Models Methods Appl. Sci.}, 25:587--645, 2015.

\bibitem{gla:hel:fle:cla:17}
Dennis Gl\"aser, Rainer Helmig, Bernd Flemisch, and Holger Class.
\newblock A discrete fracture model for two-phase flow in fractured porous
  media.
\newblock {\em Adv. Water Resour.}, 110:335--348, 2017.

\bibitem{Hackb}
Wolfgang Hackbusch.
\newblock {\em Multi-grid Methods and Applications}.
\newblock Springer, Berlin, 1985.

\bibitem{haj:kar:jen:2011}
Hadi Hajibeygi, Dimitris Karvounis, and Patrick Jenny.
\newblock A hierarchical fracture model for the iterative multiscale finite
  volume method.
\newblock {\em Journal of Computational Physics}, 230(24):8729 -- 8743, 2011.

\bibitem{hoa:jap:ker:rob:16}
Thi-Thao-Phuong Hoang, Caroline Japhet, Michel Kern, and Jean~E. Roberts.
\newblock Space-time domain decomposition for reduced fracture models in mixed
  formulation.
\newblock {\em SIAM J. Numer. Anal.}, 54(1):288--316, 2016.

\bibitem{kei:fum:ber:ste:17}
E.~Keilegavlen, A.~Fumagalli, R.~Berge, and I.~Stefansson.
\newblock Implementation of mixed-dimensional models for flow in fractured
  porous media.
\newblock {\em arXiv:1712.07392 [cs.CE]}, 2017.

\bibitem{kna:rob:14}
Peter Knabner and Jean~E. Roberts.
\newblock Mathematical analysis of a discrete fracture model coupling {D}arcy
  flow in the matrix with {D}arcy--{F}orchheimer flow in the fracture.
\newblock {\em ESAIM Math. Model. Numer. Anal.}, 48(5):1451--1472, 2014.

\bibitem{les:ang:qua:11}
Matteo Lesinigo, Carlo D'Angelo, and Alfio Quarteroni.
\newblock A multiscale {D}arcy--{B}rinkman model for fluid flow in fractured
  porous media.
\newblock {\em Numer. Math.}, 117(4):717--752, 2011.

\bibitem{PeiyaoJCP}
P.~Luo, C.~Rodrigo, F.~J. Gaspar, and C.~W. Oosterlee.
\newblock Monolithic multigrid method for the coupled {S}tokes flow and
  deformable porous medium system.
\newblock {\em J. Comput. Phys.}, 353:148--168, 2018.

\bibitem{PeiyaoSISC}
Peiyao Luo, Carmen Rodrigo, Francisco~J. Gaspar, and Cornelis~W. Oosterlee.
\newblock Uzawa smoother in multigrid for the coupled porous medium and
  {S}tokes flow system.
\newblock {\em SIAM J. Sci. Comp.}, 39(5):S633--S661, 2017.

\bibitem{Martin-Jaffre-Roberts}
Vincent Martin, Jérôme Jaffr\'e, and Jean~E. Roberts.
\newblock Modeling fractures and barriers as interfaces for flow in porous
  media.
\newblock {\em SIAM J. Sci. Comp.}, 26(5):1667--1691, 2005.

\bibitem{nor:boo:fum:kei:18}
J.~M. Nordbotten, W.~M. Boon, A.~Fumagalli, and E.~Keilegavlen.
\newblock Unified approach to discretization of flow in fractured porous media.
\newblock {\em arXiv:1802.05961 [math.NA]}, 2018.

\bibitem{LFA_overlapping}
C.~Rodrigo, F.~J. Gaspar, and F.~J. Lisbona.
\newblock On a local {F}ourier analysis for overlapping block smoothers on
  triangular grids.
\newblock {\em Appl. Numer. Math.}, 105:96--111, 2016.

\bibitem{myBook}
Carmen Rodrigo.
\newblock {\em Geometric Multigrid Methods on Triangular Grids. Application to
  Semi-structured Meshes}.
\newblock LAP Lambert Academic Publishing, Germany, 2012.

\bibitem{ARCME}
Carmen Rodrigo, Francisco~J. Gaspar, and Francisco~J. Lisbona.
\newblock Multigrid methods on semi-structured grids.
\newblock {\em Arch. Comput. Methods Eng.}, 19(4):499--538, Dec 2012.

\bibitem{rus:whe:83}
T.~F. Russell and M.~F. Wheeler.
\newblock Finite element and finite difference methods for continuous flows in
  porous media.
\newblock In R.~E. Ewing, editor, {\em The Mathematics of Reservoir
  Simulation}, volume~1 of {\em Frontiers in Applied Mathematics}, pages
  35--106. SIAM, Philadelphia, 1983.

\bibitem{san:ber:nor:12}
T.~H. Sandve, I.~Berre, and J.~M. Nordbotten.
\newblock An efficient multi-point flux approximation method for discrete
  fracture-matrix simulations.
\newblock {\em J. Comput. Phys.}, 231(9):3784--3800, 2012.

\bibitem{san:kei:nor:14}
T.~H. Sandve, E.~Keilegavlen, and J.~M. Nordbotten.
\newblock Physics-based preconditioners for flow in fractured porous media.
\newblock {\em Water Resour. Res.}, 50:1357--1373, 2014.

\bibitem{sch:fle:hel:woh:15}
Nicolas Schwenck, Bernd Flemisch, Rainer Helmig, and Barbara~I. Wohlmuth.
\newblock Dimensionally reduced flow models in fractured porous media:
  crossings and boundaries.
\newblock {\em Comput. Geosci.}, 19(6):1219--1230, 2015.

\bibitem{Stu_Tro}
Klaus St\"uben and Ulrich Trottenberg.
\newblock Multigrid methods: fundamental algorithms, model problem analysis and
  applications.
\newblock In W.~Hackbusch and U.~Trottenberg, editors, {\em Multigrid Methods},
  volume 960 of {\em Lecture Notes in Math.}, pages 1--176. Springer, Berlin,
  1982.

\bibitem{ten:saa:haj:2016}
Matei Tene, Mohammed Saad~Al Kobaisi, and Hadi Hajibeygi.
\newblock Algebraic multiscale method for flow in heterogeneous porous media
  with embedded discrete fractures (f-ams).
\newblock {\em Journal of Computational Physics}, 321:819 -- 845, 2016.

\bibitem{TOS01}
Ulrich Trottenberg, Cornelis~W. Oosterlee, and Anton Sch\"uller.
\newblock {\em Multigrid}.
\newblock Academic Press, New York, 2001.

\bibitem{vanka}
S.~P. Vanka.
\newblock Block-implicit multigrid solution of {N}avier-{S}tokes equations in
  primitive variables.
\newblock {\em J. Comput. Phys.}, 65(1):138--158, 1986.

\bibitem{Wess}
P.~Wesseling.
\newblock {\em An Introduction to Multigrid Methods}.
\newblock John Wiley \& Sons, Ltd., Chichester, 1992.

\bibitem{Wie01}
R.~Wienands and W.~Joppich.
\newblock {\em Practical Fourier Analysis for Multigrid Methods}.
\newblock Chapman \& Hall/CRC Press, Boca Raton, 2005.

\end{thebibliography}
\end{document}